\documentclass[12pt]{amsart}

\usepackage{amsmath,amssymb,amsthm,amscd,a4wide,upref}

\usepackage[enableskew]{youngtab}

\usepackage[colorlinks=true, linkcolor=blue, citecolor=magenta, urlcolor=blue, backref=page]{hyperref} 

\hypersetup{colorlinks,citecolor=blue,filecolor=black,linkcolor=blue,urlcolor=blue}

\usepackage{enumerate}

\usepackage{mathtools}

\usepackage{enumitem}

\usepackage{lmodern}     

\usepackage[T1]{fontenc}

\newcommand{\End}{{\rm{End}\ts}}

\newcommand{\Mat}{{\rm{Mat}}}

\newcommand{\sh}{{\rm{sh}}}

\newcommand{\diag}{ {\rm diag}}

\newcommand{\ot}{\otimes}

\newcommand{\la}{\lambda}

\newcommand{\La}{\Lambda}

\newcommand{\al}{\alpha}

\newcommand{\Ga}{\Gamma}

\newcommand{\si}{\sigma}

\newcommand{\vp}{\varphi}

\newcommand{\de}{\delta}

\newcommand{\ts}{\,}

\newcommand{\tss}{\hspace{1pt}}

\newcommand{\U}{ {\rm U}}

\newcommand{\CC}{\mathbb{C}\tss}

\newcommand{\ZZ}{\mathbb{Z}\tss}

\newcommand{\Ac}{\mathcal{A}}

\newcommand{\Pc}{\mathcal{P}}

\newcommand{\Qc}{\mathcal{Q}}

\newcommand{\Bc}{\mathcal{B}}

\newcommand{\Dc}{\mathcal{D}}

\newcommand{\Hc}{\mathcal{H}}

\newcommand{\Uc}{\mathcal{U}}

\newcommand{\gl}{\mathfrak{gl}}

\newcommand{\HC}{ {\mathcal HC}}

\newcommand{\tr}{ {\rm tr}}

\newcommand{\str}{ {\rm str}}

\newcommand{\bi}{\bar{\imath}}

\newcommand{\Sym}{\mathfrak S}

\newcommand{\Fand}{\qquad\text{and}\qquad}

\newcommand{\Ucirc}{U_q^\circ(\gl_{M|N})}

\newtheorem{thm}{Theorem}[section]

\newtheorem{lem}[thm]{Lemma}

\newtheorem{prop}[thm]{Proposition}

\newtheorem{cor}[thm]{Corollary}

\newtheorem{conj}[thm]{Conjecture}

\theoremstyle{definition}

\newtheorem{defin}[thm]{Definition}

\theoremstyle{remark}

\newtheorem{remark}[thm]{Remark}

\newtheorem{example}[thm]{Example}

\newcommand{\bth}{\begin{thm}}

\renewcommand{\eth}{\end{thm}}

\newcommand{\bpr}{\begin{prop}}

\newcommand{\epr}{\end{prop}}

\newcommand{\ble}{\begin{lem}}

\newcommand{\ele}{\end{lem}}

\newcommand{\bco}{\begin{cor}}

\newcommand{\eco}{\end{cor}}

\newcommand{\bde}{\begin{defin}}

\newcommand{\ede}{\end{defin}}

\newcommand{\bex}{\begin{example}}

\newcommand{\eex}{\end{example}}

\newcommand{\bre}{\begin{remark}}

\newcommand{\ere}{\end{remark}}

\newcommand{\bcj}{\begin{conj}}

\newcommand{\ecj}{\end{conj}}

\newcommand{\bal}{\begin{aligned}}

\newcommand{\eal}{\end{aligned}}

\newcommand{\beq}{\begin{equation}}

\newcommand{\eeq}{\end{equation}}

\newcommand{\ben}{\begin{equation*}}

\newcommand{\een}{\end{equation*}}

\newcommand{\bpf}{\begin{proof}}

\newcommand{\epf}{\end{proof}}

\newcommand{\BZ}{\mathbb{Z}}            

\newcommand{\Rc}{\mathcal{R}}

\def\beql#1{\begin{equation}\label{#1}}

\newcommand{\Rcheck}{\check{R}}

\title{Super $q$-Immanants and the Center of the Reflection Equation Algebra}
\author{Zheng Li}
\address{School of Artificial Intelligence, Jianghan University, Wuhan, Hubei 430056, China}
\email{lz1994@jhun.edu.cn}
\author{Ming Liu}
\address{School of Artificial Intelligence, Jianghan University, Wuhan, Hubei 430056, China}
\email{ming.l1984@gmail.com}
\author{Jian Zhang}
\address{School of Mathematics and Statistics,  Central China Normal University, Wuhan, Hubei 430079, China}
\email{jzhang@ccnu.edu.cn}
\date{}

\begin{document}

\begin{abstract}
We construct super $q$-immanants for the quantum general linear superalgebra $U_q(\gl_{M|N})$, associated with primitive idempotents of the Hecke algebra and parameterized by Young diagrams.
Their coefficients are central elements whose Harish--Chandra images are expressed in terms of factorial supersymmetric Schur polynomials. We also determine the center of the reflection
equation subalgebra generated by the entries of $L=L^+(L^-)^{-1}$.
We prove that the Harish-Chandra homomorphism identifies this center with the algebra of supersymmetric polynomials, and that the central elements associated with hook Young diagrams
form a basis.
The proof relies on a PBW theorem for the modified reflection equation algebra and a comparison with the classical center of $U(\gl_{M|N})$.
Finally, we establish higher quantum super Capelli identities through a realization of the reflection equation algebra in a braided Weyl superalgebra.
\end{abstract}

\maketitle

\section{Introduction}
\label{sec:introduction}

The quantum immanants for the general linear Lie algebra $\gl_n$ were introduced by Okounkov as a distinguished family of central elements of $U(\gl_n)$ associated with Young diagrams. 
These elements form a basis of the center, and their images under the Harish--Chandra isomorphism are
the shifted Schur polynomials, which can be regarded as special cases of factorial Schur polynomials.
Their realizations as polynomial differential operators also lead to higher Capelli identities.
The corresponding theory for the general linear Lie superalgebra $\gl_{M|N}$ was developed using factorial supersymmetric Schur polynomials and super Capelli identities; see \cite{Molev 1998}.

Quantum analogues of these constructions have attracted considerable attention.
In \cite{JLM 2024}, Jing, Liu and Molev constructed explicit quantum Sugawara operators for the quantum affine algebra $U_q(\widehat{\gl}_n)$ at the critical level, associated with primitive idempotents of the Hecke algebra. They also determined the Harish-Chandra images of these operators.
Subsequently, in \cite{JLM 2025}, the same authors introduced $q$-immanant polynomials for the quantum
group $U_q(\gl_n)$ and showed that their Harish--Chandra images are factorial Schur polynomials.
They further established higher quantum Capelli identities by means of a braided Weyl algebra.
More recently, Jing, Liu and Zhang \cite{JLZ 2026} extended the construction of quantum Sugawara operators to the quantum affine superalgebra $U_q(\widehat{\gl}_{M|N})$ and calculated their
Harish-Chandra images.

The purpose of the present paper is to study the corresponding constructions for the quantum
general linear superalgebra $U_q(\gl_{M|N})$. Besides extending the $q$-immanant construction
to the super case, we investigate the center of the reflection equation subalgebra generated by
the entries of the matrix
\[
L=L^+(L^-)^{-1}.
\]
This problem is closely related to the Harish-Chandra theory for quantum superalgebras, but requires a separate treatment.
Indeed, the reflection equation algebra considered here does not contain the inverses of its Cartan
generators, and its center need not coincide with the center of the ambient quantum superalgebra.
Reflection equation algebras of $GL(M|N)$ type and their modified versions have been studied
in \cite{GPS 2009}.
Generalized Harish--Chandra morphisms for reflection equation algebras were investigated in
\cite{GS 2025}, while the Harish-Chandra homomorphism for quantum superalgebras was studied
in \cite{LWW 2022}.
Our approach gives an explicit description of the center of the reflection equation
subalgebra through the actual PBW Cartan projection.

Throughout the paper, let $q$ be an indeterminate.
We work over $\CC(q)$ and set
\[
\Ucirc=\CC(q)\langle L_{ij}:1\leq i,j\leq M+N\rangle\subset U_q(\gl_{M|N}).
\]
Our first main result determines the center of $\Ucirc$. For each standard tableau $\La$ of shape
$\la\vdash r$, let $E_\La$ be the corresponding primitive Hecke idempotent acting on the tensor
power of the natural representation.
We construct central elements
\[
C_\La=\operatorname{str}_{q;1,\ldots,r}\left(L_{\bar1}\cdots L_{\bar r}E_\La\right),
\]
where $L_{\bar a}$ are the braided copies of $L$.
We show that $C_\La$ depends only on the shape $\la$, and write it as $C_\la$.
Moreover, $C_\la=0$ if $\la$ lies outside the $(M,N)$-hook.

Let 
\[
H_{M,N}=\{\la:\la_{M+1}\leq N\}.
\]
We prove that
\[
Z(\Ucirc)=\bigoplus_{\la\in H_{M,N}}\CC(q)C_\la,
\]
and that the PBW Cartan projection induces an algebra isomorphism
\[
\mathrm{HC}:Z(\Ucirc)\xrightarrow{\ \sim\ }\Lambda_{M|N}, \qquad
\mathrm{HC}(C_\la)=s_\la(u-v),
\]
where $\Lambda_{M|N}$ is the algebra of supersymmetric polynomials and $s_\la(u-v)$ is the supersymmetric Schur polynomial.
In particular, every central element of $\Ucirc$ is central in the ambient quantum superalgebra, so that
\[
Z(\Ucirc)=\Ucirc\cap Z(U_q(\gl_{M|N})).
\]
We also obtain a filtered version of this description, with the basis indexed by hook partitions of bounded size.

The proof of the center theorem uses the modified reflection generators
\[
X=\frac{I-L}{q-q^{-1}}.
\]
We establish a PBW theorem for an integral form of the modified reflection equation algebra, whose specialization at $q=1$ is isomorphic to $U(\gl_{M|N})$.
This allows us to compare the dimensions of the filtered centers at generic $q$ and at $q=1$.
The classical Harish-Chandra theorem, together with the hook Schur basis theorem, then gives an upper bound for the dimension of each filtered component of $Z(\Ucirc)$.
On the other hand, the Harish-Chandra images of the elements $C_\La$ are calculated through weighted supertraces of Hecke idempotents.
The super Frobenius character formula identifies these images with the supersymmetric Schur polynomials, providing the matching lower bound and completing the description of the center.

Our second main result concerns the super $q$-immanants.
For a standard tableau $\La$ of shape $\la\vdash r$, with contents
$c_1(\La),\ldots,c_r(\La)$, we define
\[
S_\La(z)= \operatorname{str}_{q;1,\ldots,r}\left(\prod_{a=1}^{r}\left(L_{\bar a}+zq^{2c_a(\La)}\right)
E_\La\right).
\]
We show that all coefficients of $S_\La(z)$ belong to the center of $U_q(\gl_{M|N})$, and that this polynomial depends only on the shape $\la$.
For a highest-weight module with highest weight
\[
(\eta_1,\ldots,\eta_M,\theta_1,\ldots,\theta_N),
\]
the eigenvalue of $S_\la(z)$ is given by the factorial supersymmetric Schur polynomial
\[
s_\la^{M|N}(x/y\mid a),
\]
where
\[
x_i=q^{2(\eta_i-i+1)},\qquad y_j=-q^{-2(\theta_j+M-j)},\qquad a_k=zq^{2(k-M+N)}.
\]
The calculation uses the Hecke-algebra character method developed in \cite{JLM 2024,JLM 2025,JLZ 2026}, together with the factorial supersymmetric Schur functions of \cite{Molev 1998}.
At $z=0$, the polynomials $S_\la(z)$ recover the homogeneous central elements $C_\la$ constructed above.

Finally, we consider higher quantum super Capelli identities.
Following the braided differential-operator approach of \cite{JLM 2025}, we introduce a braided Weyl superalgebra and construct a homomorphism from $\Ucirc$ into this algebra using the modified reflection generators.
We then obtain a normal-ordered matrix identity associated with each standard Young tableau, together with the corresponding quantum-supertrace identity.
The matrix Capelli identity for general Hecke symmetries has been investigated in \cite{GPS 2022,GSZ 2025,Zaitsev 2024}.
Our treatment specializes this framework to the standard quantum general linear superalgebra and relates it to the reflection subalgebra considered here.
We keep the Capelli normalization distinct from that of the super $q$-immanants and recover the corresponding higher super Capelli identity in the classical limit $q\to1$.

The paper is organized as follows. In Section~2, we recall the $R$-matrix realization of $U_q(\gl_{M|N})$, its PBW decomposition and the Cartan projection.
Section~3 collects the necessary results on the Hecke algebra, including primitive idempotents and the fusion procedure.
In Section~4, we recall supersymmetric Schur polynomials and the classical Harish-Chandra theory.
Section~5 is devoted to the reflection equation algebra and the determination
of its center.
In Section~6, we construct the super $q$-immanants and calculate their factorial Harish-Chandra images.
Finally, Section~7 establishes higher quantum super Capelli identities through the braided Weyl superalgebra.

\section{The $R$-matrix and the quantum superalgebra}\label{sec:R quantum}
We equip $I=\{1,\ldots,M+N\}$ with the standard parity sequence
\beql{eq:parity}
\bar i=
\begin{cases}
0,&1\leq i\leq M,\\
1,&M<i\leq M+N,
\end{cases}
\qquad q_i=q^{1-2\bar i}.
\eeq
For the structural and center arguments we take $q$ to be an indeterminate and work over $\CC(q)$; the identities in the following sections also admit the usual specializations whenever their denominators are nonzero. 
The center classification, which uses specialization at $q=1$, is asserted only for generic $q$.  
All tensor products are understood in the graded sense. 
In particular, the $R$-matrices, the $L$-operators and the Hecke operators appearing below are even.

\subsection{The $R$-matrix}
Let $E_{ij}\in\End\CC^{M|N}$ be the standard matrix units.
We use the $R$-matrix convention of \cite{JLZ 2026}. Define
\beq\label{eq:Rzw}
\begin{aligned}
R(z,w)={}&\sum_{i\in I}(zq_i-wq_i^{-1})E_{ii}\ot E_{ii}+(z-w)\sum_{i\neq j}E_{ii}\ot E_{jj}\\
&\quad+z\sum_{i<j}(q_j-q_j^{-1})E_{ij}\ot E_{ji}+w\sum_{i>j}(q_j-q_j^{-1})E_{ij}\ot E_{ji}.
\end{aligned}
\eeq
It satisfies the Yang--Baxter equation
\beql{eq:YBE}
R_{12}(z_1,z_2)R_{13}(z_1,z_3)R_{23}(z_2,z_3)=R_{23}(z_2,z_3)R_{13}(z_1,z_3)R_{12}(z_1,z_2).
\eeq
Introduce the graded permutation operator
\beql{eq:P}
P=\sum_{i,j\in I}(-1)^{\bar j}E_{ij}\ot E_{ji}.
\eeq
Set
\[
R=R(1,0),\qquad R'=PR^{-1}P.
\]
Then
\beql{eq:R-linear}
R(z,w)=zR-wR',\qquad R-R'=(q-q^{-1})P.
\eeq
For a matrix $A=\sum a_{ij}E_{ij}$ we use the supertranspose
\beql{eq:st}
A^{st}=\sum_{i,j}(-1)^{\bar i(\bar i+\bar j)}a_{ji}E_{ij}.
\eeq
Let
\beql{eq:D}
D=\diag(d_1,\ldots,d_{M+N}),\qquad
d_i=
\begin{cases}
q^{-2i+2},&1\leq i\leq M,\\
q^{-2M+2j},&i=M+j,\quad 1\leq j\leq N.
\end{cases}
\eeq
A direct calculation gives
\beql{eq:constant-crossing}
\begin{aligned}
(R_{12}^{-1})^{st_2}D_2R_{12}^{st_2}&=D_2,\\
R_{12}^{st_1}D_1(R_{12}^{-1})^{st_1}&=D_1.
\end{aligned}
\eeq

\subsection{The quantum superalgebra $\U_q(\gl_{M|N})$}
We now recall the $RLL$ realization of the quantum general linear superalgebra. 
Let $\U_q(\gl_{M|N})$ be the unital associative superalgebra generated by the entries of the two even matrices
\[
L^\pm=\sum_{i,j\in I}E_{ij}\ot l^\pm_{ij},
\]
subject to
\beql{eq:triangular}
l^+_{ji}=0,\qquad l^-_{ij}=0\quad (1\leq i<j\leq M+N),\qquad
l^+_{ii}l^-_{ii}=l^-_{ii}l^+_{ii}=1,
\eeq
and the $RLL$ relations
\beql{RLL-finite}
RL^\pm_1L^\pm_2=L^\pm_2L^\pm_1R, \qquad RL^+_1L^-_2=L^-_2L^+_1R.
\eeq
The parity of $l^\pm_{ij}$ is $\bar i+\bar j$.  

We set
\beql{eq:L}
\Rcheck=PR.
\eeq
Then in terms of $\Rcheck$, the $RLL$ relation can be written as
\beql{RcLL-finite}
\Rcheck L^\pm_1L^\pm_2=L^\pm_1L^\pm_2\Rcheck, \qquad \Rcheck L^+_1L^-_2=L^-_1L^+_2\Rcheck.
\eeq

The $RLL$ relations also imply
\beql{ev RTT}
\begin{aligned}
&R(z,w)(L_1^++z^{-1}L_1^-)(L_2^++w^{-1}L_2^-)\\
&\quad=(L_2^++w^{-1}L_2^-)(L_1^++z^{-1}L_1^-)R(z,w).
\end{aligned}
\eeq
Indeed, this follows by substituting $R(z,w)=zR-wR'$ and using the $RLL$ relations together with
\[
R-R'=(q-q^{-1})P.
\]

Define the quantum supertrace by
\beql{eq:qsupertrace}
\operatorname{str}_q(X)=\operatorname{str}(DX).
\eeq
For several tensor copies, we write
\[
\operatorname{str}_{q;1,\ldots,r}(X)=\operatorname{str}_{1,\ldots,r}(D_1\cdots D_rX).
\]

\subsection{PBW decomposition and the Cartan projection}\label{subsec:hc-pbw-prelim}
We recall that the $RLL$ realization of $U_q(\mathfrak{gl}_{M|N})$ is isomorphic to its Drinfeld--Jimbo realization \cite{Zhang 2016}. 
Let $\epsilon_1,\ldots,\epsilon_{M+N}$ be the standard weight basis, and set $\alpha_i=\epsilon_i-\epsilon_{i+1}$ for $1\leq i<M+N$.

In the Drinfeld--Jimbo realization, we denote by $E_i$, $F_i$ the Chevalley generators associated with the simple roots $\alpha_i$, $-\alpha_i$, respectively, and by $K_a^{\pm1}$ ($1\leq a\leq M+N$) the Cartan generators.
Under the identification with the $RLL$ realization, the Cartan generators correspond to
\[
K_a=l_{aa}^+=(l_{aa}^-)^{-1}.
\]
In particular, the Cartan element associated with the simple root $\alpha_i$ is
\[
K_iK_{i+1}^{-1}=l_{ii}^+l_{i+1,i+1}^-.
\]
Set
\[
\mathcal U^0=\CC(q)[K_1^{\pm1},\ldots,K_{M+N}^{\pm1}].
\]
Let $\mathcal U^+$ and $\mathcal U^-$ denote the unital subalgebras generated by the Chevalley generators $E_i$ and $F_i$, respectively. 
Equivalently, they are generated by the positive root vectors $E_\beta$ and the negative root vectors $F_\beta$, where $\beta$ ranges over all positive roots.
Under the above identification, the entries of $L^+$ correspond to negative root directions, whereas the entries of $L^-$ correspond to positive root directions, up to Cartan factors and normalization constants.

By the standard triangular decomposition and PBW theorem for quantum general linear superalgebras, multiplication induces an isomorphism of vector spaces
\beq
\mathcal U^-\otimes\mathcal U^0\otimes\mathcal U^+\xrightarrow{\ \sim\ }U_q(\mathfrak{gl}_{M|N}).
\eeq
In particular, the ordered monomials
\[
F^{\mathbf a}(l^+)^{\mathbf b}E^{\mathbf c},
\]
form a PBW basis, where $(l^+)^{\mathbf b}=\prod_{i=1}^{M+N}(l_{ii}^+)^{b_i}$ for $\mathbf b=(b_1,\ldots,b_{M+N})\in\mathbb Z^{M+N}$ and the exponents of the even root vectors are nonnegative integers,
while those of the odd root vectors belong to $\{0,1\}$.

Let $\Uc^+_+$ and $\Uc^-_+$ be spanned by the PBW monomials of positive length in the corresponding
root vectors, thus
\[
\Uc^\pm=\CC(q)\oplus\Uc^\pm_+.
\]
Consequently, the triangular decomposition gives the following direct sum decomposition of vector spaces:
\beq
\U_q(\gl_{M|N})=\Uc^0\oplus\Uc^-_+\Uc^0\oplus\Uc^0\Uc^+_+\oplus\Uc^-_+\Uc^0\Uc^+_+.
\eeq

We define the \emph{PBW Cartan projection}
\beq\label{eq:center-Cartan-projection}
\mathrm{pr}:\U_q(\gl_{M|N})\longrightarrow\Uc^0
\eeq
to be the projection onto $\Uc^0$ along the above direct sum decomposition.

The restriction of $\mathrm{pr}$ to the center $Z(U_q(\gl_{M|N}))$ is an algebra homomorphism, called the Harish--Chandra homomorphism. 

We also consider the topologically free $\hbar$-adic quantum superalgebra $\mathcal U_\hbar$, following \cite{Yamane 1994,Zhang 1998}, with
\[
\mathcal U_\hbar/\hbar\mathcal U_\hbar\simeq U(\gl_{M|N}).
\]
Under the identification of the RLL and Drinfeld--Jimbo realizations with $q=e^\hbar$, the Cartan generators and the universal $R$-matrix have the expansions
\beql{eq:formal-classical-expansion}
\begin{aligned}
l_{ii}^+&=1+\hbar(-1)^{\bar i}E_{ii}+O(\hbar^2),\\
\Rc&=1+\hbar\left(\sum_i(-1)^{\bar i}E_{ii}\ot E_{ii}+2\sum_{i<j}(-1)^{\bar j}E_{ij}\ot E_{ji}\right)+O(\hbar^2).
\end{aligned}
\eeq
Here the coefficients of $\hbar$ are identified with elements of $U(\gl_{M|N})$ through the classical limit.

The Cartan factor of $\Rc$ is
\[
\mathcal Q=\exp\left(\hbar\sum_i(-1)^{\bar i}E_{ii}\ot E_{ii}\right).
\]
Let
\[
\pi:\mathcal U_\hbar\longrightarrow\End_{\CC[[\hbar]]}(\CC^{M|N}[[\hbar]])
\]
be the natural representation. With our tensor conventions,
\beql{eq:universal-L}
L^+=(\pi\ot\mathrm{id})(\Rc),\qquad (L^-)^{-1}=(\pi\ot\mathrm{id})(\Rc_{21}).
\eeq

The substitution $q=e^\hbar$ induces injective homomorphisms $\iota$
\beql{eq:rational-to-formal}
U_q(\gl_{M|N})\hookrightarrow\mathcal U_\hbar\ot_{\CC[[\hbar]]}\CC((\hbar)),
\eeq
which sends $l_{ii}^+$ to $\exp\left(\hbar(-1)^{\bi}E_{ii}\right)$.
The injectivity follows from the PBW theorems for the rational and formal quantum superalgebras \cite{ZC 2018,Yamane 1994}; see also \cite[Section~3.3, Proposition~4(ii)]{Zhang 1998} for a related argument using the separation of elements by tensor powers of the natural representation.

Let $\varepsilon$ denote the counit of $U_q(\gl_{M|N})$, whose restrictions to $\mathcal U^\pm$ are the augmentation maps determined by
\[
\varepsilon(1)=1,\qquad \varepsilon(E_i)=\varepsilon(F_i)=0.
\]
The embedding $\iota$ preserves the positive and negative root subalgebras and their augmentations. Hence, for $F\in\mathcal U^-$, $a\in\mathcal U^0$, and $E\in\mathcal U^+$, we have
\[
\begin{aligned}
\operatorname{pr}_\hbar(\iota(FaE))&=\varepsilon(F)\iota_0(a)\varepsilon(E)\\
&=\iota_0(\operatorname{pr}(FaE)),
\end{aligned}
\]
where $\iota_0=\iota|_{\mathcal U^0}$.
By the PBW decomposition and linearity,
\[
\operatorname{pr}_\hbar\circ\iota=\iota_0\circ\operatorname{pr}.
\]

\section{The Hecke algebra and its representations}\label{sec:Hecke-algebra-prelim}
We recall the type-$A$ Hecke algebra and the primitive idempotents used in both the reflection-algebra center calculation and the $q$-immanant construction. Let $\Hc_r$ be the Hecke algebra of type $A_{r-1}$ generated by $T_1,\ldots,T_{r-1}$ with defining relations
\beq\label{eq:Hecke}
(T_i-q)(T_i+q^{-1})=0, \qquad T_iT_{i+1}T_i=T_{i+1}T_iT_{i+1}, \qquad
T_iT_j=T_jT_i\quad (|i-j|>1).
\eeq
We use the standard tensor representation
\beql{Hecke-action}
\varrho_r:\Hc_r\longrightarrow \End\bigl((\CC(q)^{M|N})^{\ot r}\bigr),\qquad T_i\longmapsto\Rcheck_i.
\eeq
The braid and quadratic relations for $\Rcheck$ ensure that this is well defined.  
For specialization at $q=1$, we use its integral version over the local ring at $q=1$, introduced in Section~\ref{subsec:center-idempotents-prelim}.

\subsection{Young basis}
Let $\lambda = (\lambda_1,\dots,\lambda_\ell)$ be a partition of $r$, identified with its Young diagram. 
Its conjugate partition is denoted $\lambda'=(\lambda_1',\lambda_2',\ldots)$, where $\lambda_j'$ counts the number of boxes in column $j$. 
For a box $\alpha = (i,j)\in \lambda$, the hook length is $h(\alpha) = \lambda_i + \lambda_j' - i - j + 1$, and the content is $c(\alpha) = j - i$.
	
A $\lambda$-tableau $\La$ is a filling of the diagram $\lambda$ with numbers taken from $\{1,\dots,r\}$. 
It is called \emph{standard} if the entries increase strictly along each row and down each column; it is \emph{semistandard} if the entries increase weakly along rows and strictly down columns.
Given a standard $\lambda$-tableau $\Lambda$, write $c_k(\Lambda)$ for the content of the cell containing the entry $k$, and define $d_k(\Lambda) = c_{k+1}(\Lambda) - c_k(\Lambda)$.
	
The irreducible representations of $\mathcal{H}_r$ over $\CC(q)$ are labelled by partitions of $r$. For a partition $\lambda \vdash r$, let $V_\lambda$ be the corresponding irreducible module and let
\beql{phila}
\vp_{\la}:\Hc_{r}\to \End V_{\la}
\eeq
be the associated algebra homomorphism. 
The space $V_\lambda$ carries an $\mathcal{H}_r$-invariant inner product $\langle \cdot,\cdot\rangle = \langle \cdot,\cdot\rangle_{V_\lambda}$ satisfying $\langle h v_1, v_2\rangle = \langle v_1, h^* v_2\rangle$ for all $h\in\mathcal{H}_r$ and $v_1,v_2\in V_\lambda$ (see \cite{DJ 1987,Hoefsmit 1974}).
	
There exists an orthonormal \emph{Young basis} $\{v_\Lambda\}$ of $V_\lambda$, indexed by the standard $\lambda$-tableaux $\Lambda$. Denote $f_\lambda = \dim V_\lambda$, which equals the number of such standard tableaux. Using the seminormal form developed in \cite{DJ 1987} and \cite{Hoefsmit 1974}, the action of the generators $T_k$ on this basis can be made completely explicit. After a suitable normalization (possibly over a finite extension
of $\CC(q)$ containing the square roots in the seminormal coefficients), for any $k \in \{1,\dots,r-1\}$ one obtains
\beql{eq:Young basis}
T_{k} \ts v^{}_{\La}=\frac{q^{d_k(\La)}}{[d_k(\La)]_q} v^{}_{\La}+\sqrt{1-\frac{1}{[d_k(\La)]_q^2}}
\ts v^{}_{\si_{k} \La},
\eeq
where $[n]_q = \frac{q^n - q^{-n}}{q - q^{-1}}$, the tableau $\sigma_k\Lambda$ results from $\Lambda$ by swapping the entries $k$ and $k+1$, and $v_{\sigma_k\Lambda}=0$ whenever $\sigma_k\Lambda$ fails to be standard.
	
For any skew diagram $\theta$ consisting of $r$ boxes, the same formulas \eqref{eq:Young basis} define a representation of $\mathcal{H}_r$ on the linear span $V_\theta$ of the vectors $v_\Lambda$ labelled by standard tableaux $\Lambda$ of shape $\theta$.
	
\subsection{Symmetrizing trace and Schur elements}
Define the {\em symmetrizing trace} on the Hecke algebra as the following linear map
\ben
\operatorname{tr}_{\Hc_r}:\Hc_{r}\to\CC(q),\qquad \operatorname{tr}_{\Hc_r}(T_{\si})=\de_{\si,e},
\een
where $e$ denotes the identity element of $\Sym_r$. This induces a {\em symmetric algebra}
structure on $\Hc_{r}$ associated with the following bilinear form
\ben
\Hc_{r}\ot\Hc_r\to\CC(q),\qquad h_1\ot h_2\mapsto\operatorname{tr}_{\Hc_r}(h_1h_2),
\een
which is symmetric and non-degenerate.  The dual basis of
$\{T_{\si}\}_{\si \in \Sym_r}$ with respect to the form is
$\{ T_{\si^{-1}}\}_{\si \in \Sym_r}$ \cite{GP 2000, Meliot 2017}.

For any $\pi\in\End V_{\la}$, define $I(\pi)\in \End V_{\la}$ by
\ben
I(\pi)=\sum_{\si\in \Sym_r}\vp_\la(T_{\si})\ts \pi\ts\vp_\la(T_{\si^{-1}}).
\een
A direct verification shows that $I(\pi)$ actually lies in $\End_{\mathcal{H}_r} V_\lambda$; therefore Schur's lemma forces it to be a scalar multiple of the identity. Consequently we obtain the key relation
\beql{iuschur}
I(\pi)=c_{\la}\tss\tr^{}_{V_{\la}}(\pi)\tss{\rm id}^{}_{V_{\la}},
\eeq
which will be used in the Hecke trace calculations. Here $c_\lambda$ is the \emph{Schur element} given by the Schur-element formula \cite{GP 2000,Meliot 2017}:
\ben
c_{\la}=\prod_{\al\in\la} q^{c(\al)}\tss [h(\al)]_q.
\een
Let $\chi_\lambda$ denote the character of the $\mathcal{H}_r$-module $V_\lambda$. Then the symmetrizing trace decomposes as
\beql{tauchi}
\operatorname{tr}_{\Hc_r}=\sum_{\la\vdash r}\frac{\chi_{\la}}{c_{\la}}.
\eeq

\subsection{Idempotents and matrix units}
The Hecke algebra $\Hc_r$ is semisimple and admits an isomorphism
\beql{isomhecke}
\Hc_r \cong \bigoplus_{\la \vdash r} \Mat_{f_{\la}}(\CC(q)).
\eeq
Via this isomorphism, the matrix units $e^{\la}_{\La,\La'} \in \Mat_{f_{\la}}(\CC(q))$ corresponding to $V_\lambda$
can be identified with elements of $\Hc_r$ by
\beql{mat-units}
e^{\la}_{\La,\La'}=\frac{1}{c_{\la}}\sum_{\si\in \Sym_r}
\langle T_{\si^{-1}}v^{}_{\La},v^{}_{\La'}\rangle^{}_{V_{\la}}T_{\si}.
\eeq
The characters satisfy the orthogonality relations: for any partitions $\lambda,\mu\vdash r$,
\beql{orthcha}
\sum_{\si\in \Sym_r}\chi_{\la}(T_{\si})\chi_{\mu}(T_{\si^{-1}})=\de_{\la,\mu}c_{\la}f_{\la}.
\eeq

The diagonal matrix units $e^{\la}_{\La}=e^{\la}_{\La \La}$ are primitive idempotents of $\Hc_r$.
Note the immediate properties of the primitive idempotents:
\ben
e^{\la}_{\La} e^{\la}_{\Ga}=0 \quad \text { if } \quad \La \neq \Ga\Fand
(e^{\la}_{\La})^{2}=e^{\la}_{\La}
\een
together with the decomposition of the identity in $\Hc_r$
\ben
1=\sum_{\la \vdash r} \sum_{\sh(\La)=\la} e^{\la}_{\La}.
\een

The following lemma given in \cite{JLM 2024} will play an important role in the proof of Theorem \ref{thm:main}.
\ble\label{lem:et}
Let $\La$ be a standard $\la$-tableau and let  $k \in\{1, \dots, r-1\}$.
Then
\ben
e^{\la}_{\La}\left(T_{k}-\frac{q^{d_k(\La)}}{[d_k(\La)]_q}\right)=\sqrt{1-\frac{1}{[d_k(\La)]_q^2}}
\ts e^{\la}_{\La,\si_{k} \La}
\een
and
\ben
\left(T_{k}+\frac{q^{-d_k(\La)}}{[d_k(\La)]_q}\right)e^{\la}_{\si_{k} \La}
=\sqrt{1-\frac{1}{[d_k(\La)]_q^2}} \ts e^{\la}_{\La,\si_{k} \La},
\een
assuming that $e^{\la}_{\La,\si_{k} \La}=0$ if the tableau $\si_{k} \La$ is not standard.
\ele

\subsection{Regular Hecke idempotents near $q=1$}\label{subsec:center-idempotents-prelim}
Put $\Ac=\CC[q,q^{-1}]_{(q-1)}$. 
The Hecke algebra $\Hc_r(\Ac)$ acts on $(\Ac^{M|N})^{\ot r}$ through $T_i\mapsto\Rcheck_i$. 
Set $J_1=1$, $J_{a+1}=T_aJ_aT_a$ and $B_a=(J_a-1)/(q-q^{-1})$. 
The identity
\[
B_a=\sum_{i<a}T_{(i,a)}\in\Hc_r(\Ac),\qquad B_a\big|_{q=1}=\sum_{i<a}(i,a)
\]
shows that no pole is present at $q=1$.  
If $\La$ is a standard tableau, let $\La^-$ be obtained by removing the box labelled $r$; if $b_\La$
is the removed box, the spectral recursion
\beql{eq:center-spectral-idempotent}
e_\La=e_{\La^-}\prod_{\substack{b\in\mathrm{Add}(\mathrm{sh}(\La^-))\\b\ne b_\La}}\frac{B_r-\vartheta(c(b))}{\vartheta(c(b_\La))-\vartheta(c(b))},\qquad
\vartheta(c)=\frac{q^{2c}-1}{q-q^{-1}},
\eeq
starts with $e_\varnothing=1$. 
The contents of the addable boxes are distinct. 
Every denominator is a unit of $\Ac$ because its value at $q=1$ is a nonzero integer. 
Consequently, the primitive idempotent $e_\La$ defined above belongs to $\Hc_r(\Ac)$ and specializes at $q=1$ to the corresponding primitive idempotent in $\mathbb C[\mathfrak S_r]$.
We put $E_\La=\varrho_r(e_\La)$.
The operator $E_\La$ is even, preserves tensor content, and commutes with the quantum-group action. 
Its regularity at $q=1$ is essential for the trace argument below.

\subsection{Fusion procedure}
For $k=1,\ldots,r-1$ introduce the baxterized elements
\beql{eq:Tkxy}
T_k(x,y)=T_k+\frac{q-q^{-1}}{xy^{-1}-1}.
\eeq
They satisfy
\beql{eq:baxter-YB}
T_k(x,y)T_{k+1}(x,z)T_k(y,z)=T_{k+1}(y,z)T_k(x,z)T_{k+1}(x,y)
\eeq
and
\beql{eq:baxter-unitarity}
T_k(x,y)T_k(y,x)=1-\frac{(q-q^{-1})^2xy}{(x-y)^2}.
\eeq

Let $\la\vdash r$ and let $\La$ be a standard $\la$-tableau.
Recall that $c_a(\La)$ is the content of the box occupied by $a$. 
Equip the set of pairs $(i,j)$ with $1\leq i<j\leq r$ with the ordering in which $(i,j)$
precedes $(i',j')$ if $j<j'$, or if $j=j'$ and $i<i'$.  
Put
\beql{eq:TLa}
T_\La(z_1,\ldots,z_r)=\overrightarrow{\prod_{(i,j)}}T_{j-i}(z_i,z_j).
\eeq
Let $w_r$ be the longest element of $\Sym_r$.
The fusion procedure gives the primitive idempotent
\beql{eq:fusion}
e_\La=\frac{1}{c_{\la'}}T_\La(z_1,\ldots,z_r)T_{w_r}^{-1}\bigg|_{z_1=q^{-2c_1(\La)}}\cdots
\bigg|_{z_r=q^{-2c_r(\La)}}.
\eeq

\section{Supersymmetric Schur functions and classical Harish--Chandra theory}\label{sec:center-super-Schur-prelim}
Let $u_1,\ldots,u_M,v_1,\ldots,v_N$ be independent commuting variables. 
The algebra of supersymmetric polynomials is
\beql{eq:center-supersymmetric-ring}
\Lambda_{M|N}=\bigl\{f\in\CC(q)[u_1,\ldots,u_M,v_1,\ldots,v_N]^{\Sym_M\times\Sym_N}:
f\big|_{u_i=v_j=t}\text{ is independent of }t\bigr\}.
\eeq
For $r\ge1$ put
\[
p_r(u-v)=\sum_{i=1}^M u_i^r-\sum_{j=1}^N v_j^r.
\]
The complete supersymmetric functions are determined by
\beql{eq:center-generating-h}
\sum_{r\ge0}h_r(u-v)t^r=\exp\left(\sum_{k\ge1}\frac{p_k(u-v)}{k}t^k\right)
=\frac{\prod_{j=1}^N(1-v_jt)}{\prod_{i=1}^M(1-u_it)}, \qquad h_0=1,\quad h_r=0\ (r<0).
\eeq
For a partition $\la$, define
\beql{eq:center-super-JT}
s_\la(u-v)=\det\bigl[h_{\la_i-i+j}(u-v)\bigr]_{1\le i,j\le\ell},\qquad \ell\ge\ell(\la).
\eeq
Set
\[
H_{M,N}=\{\la:\la_{M+1}\le N\},\qquad h_{M,N}^s=\#\{\la\in H_{M,N}:|\la|\le s\}.
\]
We recall the classical hook Schur basis theorem \cite{BereleRegev 1987,Macdonald 1992}.
\bpr\label{prop:center-classical-input}
The polynomial $s_\la(u-v)$ vanishes for $\la\notin H_{M,N}$,
and the polynomials $\{s_\la(u-v):\la\in H_{M,N}\}$ form a homogeneous
basis of $\Lambda_{M|N}$, with degree $|\la|$. Consequently,
\beql{eq:center-classical-hook-basis}
\dim_{\CC(q)}(\Lambda_{M|N})_{\le s}=h_{M,N}^s,\qquad \Lambda_{M|N}=\CC(q)[p_1(u-v),p_2(u-v),\ldots].
\eeq
\epr

The specialization $p_k\mapsto p_k(u-v)$ of ordinary symmetric
functions sends $h_r$ and $s_\la$ to $h_r(u-v)$ and $s_\la(u-v)$,
respectively. The Frobenius character formula therefore gives
\beql{eq:super-Frobenius}
s_\la(u-v)=\frac1{r!}\sum_{\sigma\in\Sym_r}
\chi^\la(\sigma^{-1})\prod_{a=1}^{t(\sigma)}p_{\ell_a(\sigma)}(u-v),
\qquad \la\vdash r,
\eeq
where $\ell_a(\sigma)$ are the cycle lengths of $\sigma$.

We also recall the classical Harish--Chandra theory for $U(\gl_{M|N})$.
Let $\mathfrak h$ be the Cartan subalgebra of diagonal matrices, and let
\[
\mathrm{pr}_{\mathrm{cl}}: U(\gl_{M|N})\longrightarrow U(\mathfrak h)
\]
be the Cartan projection associated with the classical triangular decomposition.
For a highest weight$(\eta_1,\ldots,\eta_M|\theta_1,\ldots,\theta_N)$, let $\xi$ be a corresponding highest weight vector satisfying
\[
E_{ii}\xi=\eta_i\xi,\quad 1\leq i\leq M,\qquad E_{M+j,M+j}\xi=\theta_j\xi,\quad 1\leq j\leq N.
\]
Then every central element $z\in Z(U(\gl_{M|N}))$ acts on $\xi$ by
\[
z\xi=\mathrm{pr}_{\mathrm{cl}}(z)(\eta_1,\ldots,\eta_M,\theta_1,\ldots,\theta_N)\xi.
\]

By the classical Harish--Chandra theorem \cite[Section~7]{Molev 1998}, the Cartan projection is
injective on $Z(U(\gl_{M|N}))$ and has supersymmetric polynomial image in the shifted highest-weight coordinates
\beql{eq:center-classical-shifts}
U_i=\eta_i-i+1,\qquad V_j=-\theta_j-M+j.
\eeq

For $s\geq0$, let $F_sU(\gl_{M|N})$ denote the classical PBW filtration, namely
\[
F_sU(\gl_{M|N})=\mathrm{span}_{\CC}\{X_1\cdots X_r:0\leq r\leq s,\ X_i\in\gl_{M|N}\}.
\]
The Cartan projection does not increase the PBW degree, and the shifts in \eqref{eq:center-classical-shifts} preserve the total polynomial degree bound.
Together with the injectivity of the classical Harish--Chandra homomorphism and Proposition~\ref{prop:center-classical-input}, this gives
\beql{eq:center-classical-bound}
\dim_{\CC}\left(Z(U(\gl_{M|N}))\cap F_sU(\gl_{M|N})\right)\leq h_{M,N}^s.
\eeq

\section{The reflection equation algebra and its Harish--Chandra isomorphism}
\label{sec:center-reflection-HC}

This section determines the \emph{center} of the reflection subalgebra generated by the entries of the matrix $L=L^+(L^-)^{-1}$. 
The parameter is generic, $q$ is an indeterminate, and all center statements are over $\CC(q)$.
The classical inputs and normalizations were recorded in Section~\ref{sec:R quantum}. 
Importantly, we do not identify the center of this subalgebra with the center of the entire quantum
superalgebra in advance. 
The latter inclusion will emerge from the filtered dimension argument.

\subsection{A PBW basis for the modified reflection generators}\label{subsec:mod-ref-gen}
Set
\beql{eq:center-A}
L=L^+(L^-)^{-1},\qquad
\Ucirc=\CC(q)\langle L_{ij}:1\leq i,j\leq M+N\rangle
\subset\U_q(\gl_{M|N}).
\eeq
Let us define
\[
L_{\bar1}=L_1,\qquad L_{\overline{a+1}}=\Rcheck_aL_{\bar a}\Rcheck_a^{-1}.
\]
The following lemma also explains why the ordered products used later belong to $\Ucirc$.

\ble\label{lem:center-RE-fusion}
The matrix $L$ satisfies the graded reflection equation
\beql{eq:center-reflection-equation}
\Rcheck L_1\Rcheck L_1=L_1\Rcheck L_1\Rcheck.
\eeq
For mutually commuting scalar parameters $t_1,\ldots,t_r$ one has
\beql{eq:center-shifted-product}
\prod_{a=1}^r(L_{\bar a}+t_a)=\left(\prod_{a=1}^r(L_a^++t_aL_a^-)\right)(L_r^-)^{-1}\cdots(L_1^-)^{-1},
\eeq
where The factors in both products are ordered by increasing index $a$.
\ele
\bpf
We first prove \eqref{eq:center-shifted-product} for arbitrary commuting parameters $t_1,\ldots,t_r$.
The proof proceeds by reverse induction on $k$ starting from the trivial case $k=r$ to show that
\begin{align*}
(L_k^++t_kL_k^-)\cdots &(L_r^++t_rL_r^-)(L_r^-)^{-1}\cdots (L_{k}^-)^{-1}=\\
&=(L_k+t_k)\Rcheck_k(L_k+t_{k+1})\Rcheck_k^{-1}\cdots \Rcheck_{r-1}\cdots \Rcheck_{k}(L_k+t_r)\Rcheck_{k}^{-1}\cdots \Rcheck_{r-1}^{-1}
\end{align*}
The induction step will be completed by observing that the following relation holds
\[
L^-_{k-1}L_k =\Rcheck_{k-1}L_{k-1}\Rcheck_{k-1}^{-1}L^-_{k-1},
\]
which is a result of the $RLL$ relation \eqref{RcLL-finite}.

Taking $r=2$ and $t_1=t_2=0$, we obtain
\[
L_1\Rcheck L_1\Rcheck^{-1}=L_1^+L_2^+(L_2^-)^{-1}(L_1^-)^{-1}.
\]
Multiply both sides by $\Rcheck$ and using \eqref{RcLL-finite} we get \eqref{eq:center-reflection-equation}.
\epf

Introduce
\beql{eq:center-modified-X}
X=\frac{I-L}{q-q^{-1}}.
\eeq
Inserting this into \eqref{eq:center-reflection-equation} and using $\Rcheck^2=I+(q-q^{-1})\Rcheck$ gives,
\beql{eq:center-modified-RE}
\Rcheck X_1\Rcheck X_1-X_1\Rcheck X_1\Rcheck
=\Rcheck X_1-X_1\Rcheck.
\eeq
We shall use the resulting modified reflection equation to define an integral algebra over $\Ac$, without assuming its torsion-freeness in advance.

Let $\Bc_{\Ac}$ be the unital associative $\Ac$-superalgebra with homogeneous generators $X_{ij}$, of parity $\bar i+\bar j$, and relations given by all matrix coefficients of \eqref{eq:center-modified-RE}. 
It has its word-length filtration $F_s\Bc_{\Ac}$.
\ble\label{lem:center-special-fiber}
There is a filtered superalgebra isomorphism
\beql{eq:center-special-fiber}
\Bc_{\Ac}/(q-1)\Bc_{\Ac}\overset{\sim}{\longrightarrow}
U(\gl_{M|N}),\qquad X_{ij}\longmapsto-(-1)^{\bar j}E_{ji}.
\eeq
\ele
\bpf
At $q=1$, $\Rcheck$ specializes to the graded flip $P$, and $PX_1P=X_2$.
Equation~\eqref{eq:center-modified-RE} consequently becomes
\[
X_1X_2-X_2X_1=X_1P-PX_1.
\]
Comparing the matrix coefficients, we obtain
\[
[X_{ij},X_{k\ell}]_s=\delta_{jk}(-1)^{(\bar i+\bar j)(\bar j+\bar\ell)+\bar j}
X_{i\ell}-\delta_{i\ell}(-1)^{\bar i}X_{kj},
\]
where $[a,b]_s=ab-(-1)^{|a||b|}ba$. Now set $G_{ij}=-(-1)^{\bar i}X_{ji}$. 
A direct substitution produces
\[
[G_{ij},G_{k\ell}]_{\rm s}=\delta_{jk}G_{i\ell}-(-1)^{(\bar i+\bar j)(\bar k+\bar\ell)}
\delta_{i\ell}G_{kj}.
\]
These are exactly the defining enveloping-superalgebra relations of the standard matrix units.
Since the change $X_{ij}\leftrightarrow G_{ji}$ is invertible and preserves word length, it identifies the two presentations and their filtrations.
\epf

The local integral form is the technical ingredient that permits comparison with the classical center.  In particular, the unmodified homogeneous reflection equation would specialize to a supercommutative algebra and would not provide the required filtered comparison.
\bth\label{thm:center-local-PBW}
Fix a total order on the generators $X_{ij}$. The corresponding ordered monomials, with nonnegative powers for even generators and powers $0$ or $1$ for odd generators, form an $\Ac$-basis of
$\Bc_{\Ac}$.
Those of length at most $s$ form an $\Ac$-basis of $F_s\Bc_{\Ac}$.
Moreover,
\beql{eq:center-saturation}
(q-1)\Bc_{\Ac}\cap F_s\Bc_{\Ac}=(q-1)F_s\Bc_{\Ac}.
\eeq
The rational specialization $\Bc_{\Ac}\ot_{\Ac}\CC(q)$ maps isomorphically onto $\Ucirc$ by
$X_{ij}\mapsto(\delta_{ij}-L_{ij})/{(q-q^{-1})}$. 
In particular,
\beql{eq:center-filter-specialization}
F_s\Bc_{\Ac}/(q-1)F_s\Bc_{\Ac}
\cong F_sU(\gl_{M|N})
\eeq
and every $F_s\Bc_{\Ac}$ is finite free.
\eth
\bpf
We first establish a PBW basis compatible with the word-length filtration.
This will also show that each filtered component is a direct summand of $\mathcal B_{\Ac}$, and hence that the filtration is saturated.

Let $W_s$ be the free $\Ac$-module spanned by all free words in the $X_{ij}$ of length at most $s$.  It is finite free. 
Let $\Pc_s$ be the set of proposed ordered PBW monomials of length at most $s$.
By Lemma~\ref{lem:center-special-fiber} and classical super PBW, the images of $\Pc_s$ form a basis after putting $q=1$.
Classical PBW reduction is obtained from quadratic-linear supercommutator relations without increasing word length. 
Thus every free word $w$ of length at most $s$ has a reduction
\beql{eq:center-classical-reduction}
w-\mathrm{NF}(w)=\sum_\nu a_\nu u_\nu r_\nu(1)v_\nu,\qquad |u_\nu|+|v_\nu|+2\le s,
\eeq
where the $r_\nu(1)$ are matrix-coefficient defining relations of the special fiber.
This is also an exact description of the kernel of $W_s/(q-1)W_s\to F_sU(\gl_{M|N})$: if a linear
combination maps to zero, its classical normal form is zero and \eqref{eq:center-classical-reduction} expresses it in the truncated relation space.

Choose finitely many of the elements $u_\nu r_\nu(1)v_\nu$ so that, together with $\Pc_s$, they form a basis of $W_s/(q-1)W_s$.
Replace each $r_\nu(1)$ by the coefficient $r_\nu(q)$ of \eqref{eq:center-modified-RE} and lift the
chosen vectors to $W_s$.
The resulting square change-of-basis matrix has determinant nonzero modulo $q-1$.
Since $\Ac$ is the local ring with maximal ideal $(q-1)$, the determinant is a unit.
These lifts together with $\Pc_s$ therefore form an $\Ac$-basis of $W_s$. 
The lifted relations vanish in $\Bc_{\Ac}$, so the proposed ordered monomials span $F_s\Bc_{\Ac}$.  Notice that no torsion-freeness assertion about $\Bc_{\Ac}$ has been assumed here.

It remains to prove independence.  In the formal realization
$q=e^\hbar$ we have, by \eqref{eq:formal-classical-expansion} and \eqref{eq:universal-L},
\beq
L=(\pi\ot\mathrm{id})(\Rc\Rc_{21})
=I+2\hbar\sum_{i,j}(-1)^{\bar j}E_{ij}\ot E_{ji}
+O(\hbar^2),
\label{eq:center-formal-L-expansion}
\eeq
and $q-q^{-1}=2\hbar+O(\hbar^3)$.  Hence the entries of $X=(I-L)/(q-q^{-1})$ lie in the formal algebra and satisfy
\beql{eq:center-formal-X-specialization}
X_{ij}\equiv-(-1)^{\bar j}E_{ji}\pmod\hbar.
\eeq
Suppose a nonzero finite $\Ac$-linear combination of ordered monomials maps to zero in the formal realization. 
Divide all coefficients by the smallest common power of $q-1$.
Cancellation in the target is legitimate because $e^\hbar-1$ is $\hbar$ times a unit and the formal algebra is $\hbar$-torsion-free.
At least one remaining coefficient is nonzero at $q=1$.
Reduction modulo $\hbar$ then yields a nontrivial relation among the distinct classical ordered super PBW monomials, contradicting the classical PBW theorem. 
Thus the proposed monomials are independent and the formal realization is injective.

The bases of the $F_s\Bc_{\Ac}$ are nested subsets of the same ordered-monomial basis, so each filtered piece is a direct summand.
Expanding an element in this basis shows that $(q-1)b\in F_s\Bc_{\Ac}$ forces $b\in F_s\Bc_{\Ac}$;
this proves \eqref{eq:center-saturation}.
Finally, the formal realization factors through the rational quantum algebra, so the rational map is injective.
Since $L=I-(q-q^{-1})X$, its image is exactly the unital subalgebra generated by the $L_{ij}$.
The PBW bases and Lemma~\ref{lem:center-special-fiber} establish \eqref{eq:center-filter-specialization}.
\epf

\subsection{Central elements and their Harish–Chandra images}\label{subsec:center-HC-ima}
For a standard tableau $\La$ with $r$ boxes set
\beql{eq:center-CT}
C_\La=\str_{q,1,\ldots,r}\left(L_{\bar1}\cdots L_{\bar r} E_\La\right), \qquad E_\La=\varrho_r(e_\La).
\eeq
We use $C_\La$ for this homogeneous element to avoid conflating it with the later
inhomogeneous immanant $S_\La(z)$; in fact $C_\La=S_\La(0)$.

\ble\label{lem:center-monodromy}
The elements $C_\La$ belong to $\Ucirc\cap Z(\U_q(\gl_{M|N}))$ and are even.
\ele
\bpf
We first prove the centrality assertion for arbitrary $h\in\mathcal H_r$ and commuting parameters
$t_1,\ldots,t_r$. 
Put
\[
M(t_1,\ldots,t_r)=\prod_{a=1}^r(L_{\bar a}+t_a).
\]
By \eqref{eq:center-shifted-product}, we have
\[
M(t_1,\ldots,t_r)=\left(\prod_{a=1}^r(L_a^++t_aL_a^-)\right)(L_r^-)^{-1}\cdots(L_1^-)^{-1}.
\]

Introduce an additional auxiliary space labelled by $0$. Set
\[
\widetilde R=(R_{21})^{-1},\qquad R_0^+=R_{0r}\cdots R_{01},\qquad
R_0^-=\widetilde R_{0r}\cdots\widetilde R_{01}.
\]
The $RLL$ relations imply
\beql{eq:center-external-RLL}
L_0^\epsilon M(t_1,\ldots,t_r)=(R_0^\epsilon)^{-1}M(t_1,\ldots,t_r)R_0^\epsilon L_0^\epsilon,
\qquad \epsilon=\pm.
\eeq
Indeed, for $\epsilon=+$, the same-sign and mixed $RLL$ relations give
\[
R_{0a}L_0^+(L_a^++t_aL_a^-)=(L_a^++t_aL_a^-)L_0^+R_{0a}.
\]
Iterating this relation yields
\[
L_0^+\prod_{a=1}^r(L_a^++t_aL_a^-)=(R_0^+)^{-1}\left(\prod_{a=1}^r(L_a^++t_aL_a^-)\right)L_0^+R_0^+.
\]
The same calculation with $L_a^-$ gives
\[
L_0^+(L_1^-\cdots L_r^-)=(R_0^+)^{-1}(L_1^-\cdots L_r^-)L_0^+R_0^+.
\]
Combining these identities proves \eqref{eq:center-external-RLL} for $\epsilon=+$.
For $\epsilon=-$, the corresponding calculation uses the equivalent $RLL$ relations with
$\widetilde R=(R_{21})^{-1}$.

By the Yang--Baxter equation,
\[
\Rcheck_a R_{0,a+1}R_{0a}=R_{0,a+1}R_{0a}\Rcheck_a,
\]
and the analogous relation holds for $\widetilde R$. Consequently,
\beql{eq:center-external-Hecke}
[R_0^\epsilon,\varrho_r(h)]=0,\qquad h\in\mathcal H_r,\quad\epsilon=\pm.
\eeq

We next use the following numerical quantum-supertrace identity:
\beql{eq:center-numerical-trace}
\str_{q;1}\left(R_{01}^{-1}Y_1R_{01}\right)=I_0\str_{q;1}(Y_1),\qquad \str_{q;1}\left(\widetilde R_{01}^{-1}Y_1\widetilde R_{01}\right)=I_0\str_{q;1}(Y_1),
\eeq
where $Y$ is an even matrix with coefficients in a superalgebra and does not act in auxiliary
space $0$. 
This identity can be checked on the matrix units. The diagonal coefficients give $(-1)^{\bar a}d_a I_0$. The possible off-diagonal coefficients are proportional to
\[
\frac{d_a}{q_a}-q_bd_b-(q-q^{-1})\sum_{a<j<b}(-1)^{\bar j}d_j,\qquad a<b,
\]
which vanishes by
\[
\frac{d_i}{q_i}=q_{i+1}d_{i+1},\qquad q_j-q_j^{-1}=(q-q^{-1})(-1)^{\bar j}.
\]
For $\widetilde R$ the corresponding calculation has the opposite ordering of indices and gives
the same result. 
Thus \eqref{eq:center-numerical-trace} holds.

Applying \eqref{eq:center-numerical-trace} successively in auxiliary spaces
$r,r-1,\ldots,1$, we obtain
\beql{eq:center-fused-trace}
\str_{q;1,\ldots,r}\left((R_0^\epsilon)^{-1}YR_0^\epsilon\right)=I_0\str_{q;1,\ldots,r}(Y)
\eeq
for every even $Y$ which does not act in auxiliary space $0$.

Using \eqref{eq:center-external-RLL}, \eqref{eq:center-external-Hecke}, and \eqref{eq:center-fused-trace}, we obtain
\[
\begin{aligned}
L_0^\epsilon \str_{q;1,\ldots,r}\left(M(t_1,\ldots,t_r)\varrho_r(h)\right)={}&\str_{q;1,\ldots,r}\left(L_0^\epsilon M(t_1,\ldots,t_r)\varrho_r(h)\right)\\
={}&\str_{q;1,\ldots,r}\left((R_0^\epsilon)^{-1}M(t_1,\ldots,t_r)R_0^\epsilon\varrho_r(h)\right)L_0^\epsilon\\
={}&\str_{q;1,\ldots,r}\left((R_0^\epsilon)^{-1}M(t_1,\ldots,t_r)\varrho_r(h)R_0^\epsilon
\right)L_0^\epsilon\\
={}&\str_{q;1,\ldots,r}\left(M(t_1,\ldots,t_r)\varrho_r(h)\right)L_0^\epsilon.
\end{aligned}
\]
Since the entries of $L^+$ and $L^-$ generate $U_q(\gl_{M|N})$, this proves that
\[
\str_{q;1,\ldots,r}\left(\prod_{a=1}^r(L_{\bar a}+t_a)\varrho_r(h)\right)\in Z(U_q(\gl_{M|N}))[t_1,\ldots,t_r].
\]
In particular, taking $h=e_\La$ and $t_1=\cdots=t_r=0$ gives the centrality of $C_\La$.
Moreover, $C_\La$ belongs to $\Ucirc$ by its definition, and it is even since the total matrices in the supertrace are even.
\epf

Let ${\boldsymbol a}=(a_1,\ldots,a_{M+N})\in\ZZ_{\ge0}^{M+N}$ with $|{\boldsymbol a}|=r$. 
Denote by $(V^{\ot r})_{\boldsymbol a}$ the span of pure tensors containing $a_i$ copies of the $i$-th basis vector, and by $\Pi_{\boldsymbol a}$ the corresponding projector. 
This content space is homogeneous of parity $\bar{\boldsymbol a}=\sum_{j=1}^N a_{M+j}\pmod2$.
Since $\Rcheck$ only exchanges tensor letters (possibly adding a multiple of the original tensor), every Hecke operator preserves these spaces.

\ble\label{lem:center-Cartan-content}
With the PBW Cartan projection of \eqref{eq:center-Cartan-projection},
\beql{eq:center-Cartan-monodromy}
(\mathrm{id}\ot\mathrm{pr})\left(L_{\bar1}\cdots L_{\bar r}\right)=\sum_{|\boldsymbol a|=r}\Pi_{\boldsymbol a}\ot\prod_{i=1}^{M+N}(l_{ii}^+)^{2a_i}.
\eeq
Consequently,
\beql{eq:center-image-content}
\mathrm{pr}(C_\La)=\sum_{|\boldsymbol a|=r}(-1)^{\bar{\boldsymbol a}}
\tr\left(E_\La\big|_{(V^{\ot r})_{\boldsymbol a}}\right)\prod_{i=1}^{M+N}(d_i (l_{ii}^+)^2)^{a_i}.
\eeq
\ele

\bpf
We first establish the identity
\beql{eq:center-monodromy}
L_{\bar1}\cdots L_{\bar r}=(\pi_r\ot\mathrm{id})(\Rc\Rc_{21}),\qquad
\pi_r=\pi^{\ot r}\circ\Delta^{(r-1)}.
\eeq
For $r=1$, this follows from \eqref{eq:universal-L}.

For $r=2$, by the coproduct identities
\[
(\Delta\ot\mathrm{id})(\Rc)=\Rc_{13}\Rc_{23},\qquad
(\Delta\ot\mathrm{id})(\Rc_{21})=\Rc_{32}\Rc_{31},
\]
we obtain
\[
(\Delta\ot\mathrm{id})(\Rc\Rc_{21})=\Rc_{13}\Rc_{23}\Rc_{32}\Rc_{31}.
\]
On the other hand, using $\Rcheck_{12}=P_{12}\Rc_{12}$ and the universal Yang--Baxter equation, we have
\[
\begin{aligned}
L_{\bar1}L_{\bar2}
&=L_1\Rcheck_{12}L_1\Rcheck_{12}^{-1}\\
&=\Rc_{13}\Rc_{31}\Rc_{21}\Rc_{23}\Rc_{32}\Rc_{21}^{-1}\\
&=\Rc_{13}\Rc_{23}\Rc_{21}\Rc_{31}\Rc_{32}\Rc_{21}^{-1}\\
&=\Rc_{13}\Rc_{23}\Rc_{32}\Rc_{31}.
\end{aligned}
\]
Here the first two tensor factors are evaluated in the natural representation. More generally,
iterating the coproduct identities gives
\[
(\Delta^{(r-1)}\ot\mathrm{id})(\Rc\Rc_{21})=\Rc_{1,r+1}\cdots\Rc_{r,r+1}\Rc_{r+1,r}\cdots\Rc_{r+1,1}.
\]
Repeated application of the same Yang--Baxter relations yields
\[
\begin{aligned}
(\pi_r\ot\mathrm{id})(\Rc\Rc_{21})=&\prod_{a=1}^{r}\left(\Rcheck_{a-1}\cdots\Rcheck_1
L_1\Rcheck_1^{-1}\cdots\Rcheck_{a-1}^{-1}\right)\\
=&L_{\bar1}\cdots L_{\bar r}.
\end{aligned}
\]
This proves \eqref{eq:center-monodromy}.

We now work in the formal quantum superalgebra.
Recall that the Cartan factor of $\Rc$ is
\[
\Qc=\exp\left(\hbar\sum_{i=1}^{M+N}(-1)^{\bi}E_{ii}\ot E_{ii}\right).
\]
Since $\Rc\Rc_{21}$ is in negative--Cartan--positive PBW order in its second tensor factor, we have
\[
(\mathrm{id}\ot\mathrm{pr}_{\hbar})(\Rc\Rc_{21})=\Qc^2.
\]
Consequently, by \eqref{eq:center-monodromy},
\[
\begin{aligned}
(\mathrm{id}\ot\mathrm{pr}_{\hbar})
(L_{\bar1}\cdots L_{\bar r})=&(\pi_r\ot\mathrm{id})(\Qc^2)\\
=&\exp\left(2\hbar\sum_{i=1}^{M+N}(-1)^{\bi}\pi_r(E_{ii})\ot E_{ii}\right).
\end{aligned}
\]

For each content $\boldsymbol a$, we have
\[
\pi_r(E_{ii})\Pi_{\boldsymbol a}=a_i\Pi_{\boldsymbol a}.
\]
It follows that
\[
(\mathrm{id}\ot\mathrm{pr}_{\hbar})(L_{\bar1}\cdots L_{\bar r})=\sum_{|\boldsymbol a|=r}\Pi_{\boldsymbol a}\ot\exp\left(2\hbar\sum_{i=1}^{M+N}(-1)^{\bi}a_iE_{ii}\right).
\]
Since the homomorphism \eqref{eq:rational-to-formal} sends
\[
l_{ii}^+\mapsto\exp\left(\hbar(-1)^{\bi}E_{ii}\right),
\]
the compatibility of the two Cartan projections and the injectivity of this homomorphism give
\[
(\mathrm{id}\ot\mathrm{pr})(L_{\bar1}\cdots L_{\bar r})=\sum_{|\boldsymbol a|=r}
\Pi_{\boldsymbol a}\ot\prod_{i=1}^{M+N}(l_{ii}^+)^{2a_i}.
\]
This proves \eqref{eq:center-Cartan-monodromy}.

Finally, since $E_\La$ preserves each content space $(V^{\ot r})_{\boldsymbol a}$, we obtain
\[
\begin{aligned}
\mathrm{pr}(C_\La)&=\str_{q;1,\ldots,r}\left((\mathrm{id}\ot\mathrm{pr})(L_{\bar1}\cdots L_{\bar r})E_\La\right)\\
&=\sum_{|\boldsymbol a|=r}\str_{V^{\ot r}}\left(D^{\ot r}\Pi_{\boldsymbol a}E_\La\right)
\prod_{i=1}^{M+N}(l_{ii}^+)^{2a_i}.
\end{aligned}
\]

Recall that
\[
D=\diag(d_1,\ldots,d_{M+N}),\qquad d_i=
\begin{cases}
q^{-2i+2},&1\leq i\leq M,\\
q^{-2M+2j},&i=M+j,\quad 1\leq j\leq N.
\end{cases}
\]
On $(V^{\ot r})_{\boldsymbol a}$, we have
\[
D^{\ot r}=\prod_{i=1}^{M+N}d_i^{a_i}\,\mathrm{id},\qquad
(-1)^{|\,\cdot\,|}=(-1)^{\bar{\boldsymbol a}}.
\]
Hence
\[
\str_{V^{\ot r}}\left(D^{\ot r}\Pi_{\boldsymbol a}E_\La\right)=(-1)^{\bar{\boldsymbol a}}
\left(\prod_{i=1}^{M+N}d_i^{a_i}\right)
\tr\left(E_\La\big|_{(V^{\ot r})_{\boldsymbol a}}\right).
\]
Therefore
\[
\mathrm{pr}(C_\La)=\sum_{|\boldsymbol a|=r}(-1)^{\bar{\boldsymbol a}}\tr\left(E_\La\big|_{(V^{\ot r})_{\boldsymbol a}}\right)\prod_{i=1}^{M+N}\left(d_i(l_{ii}^+)^2\right)^{a_i}.
\]
This proves \eqref{eq:center-image-content}.
\epf

Put $Z=\mathrm{diag}(u_1,\ldots,u_M,v_1,\ldots,v_N)$ for $u_i,v_j$ being algebraically independent formal variables.
For $h\in\Hc_r$ define
\beql{eq:center-Phi}
\Phi_r(h)=\str_{V^{\ot r}}\big(Z^{\ot r}\varrho_r(h)\big).
\eeq
The numerical operator $Z^{\ot r}$ commutes with the Hecke action, so $\Phi_r(ab)=\Phi_r(ba)$ by cyclicity.
Moreover, specialize the independent variables to Cartan elements by
\beql{eq:center-uv}
u_i=d_i (l_{ii}^+)^2,\qquad v_j=d_{M+j}(l_{M+j,M+j}^+)^2.
\eeq
Then \eqref{eq:center-image-content} says precisely that
\beql{eq:center-pr-phi}
\mathrm{pr}(C_\La)=\Phi_r(e_\La).
\eeq
The following regularity observation is useful; replacing an arbitrary Hecke operator by its value at $q=1$ would not be justified.
\ble\label{lem:center-rank-idempotent}
For any content space and any standard tableau $\La$,
\[
\tr\big(E_\La|_{(V^{\ot r})_{\boldsymbol a}}\big)\in\ZZ_{\ge0}
\]
is constant as an element of $\CC(q)$ and equals the rank of the corresponding specialized idempotent at $q=1$.
\ele
\bpf
The restricted idempotent is a matrix $B(q)$ with entries in $\Ac$.  Over the fraction field it satisfies $B(q)^2=B(q)$.
Its minimal polynomial divides $t(t-1)$, which has distinct roots in characteristic zero; hence it is diagonalizable with eigenvalues $0$ and $1$.
Its trace equals its rank, an integer $r_0$.
Taking $q=1$ preserves idempotence because all entries are regular, and preserves the trace, since the trace is the constant scalar $r_0$.
The specialized rank equals that same integer.
\epf

\bpr\label{prop:center-weighted-Frob}
For any standard tableau $\La$ of shape $\la\vdash r$,
\beq\label{eq:center-Phi-Schur}
\Phi_r(e_\La)=s_\la(u-v).
\eeq
More generally, if $\chi_q^\la$ is the irreducible Hecke character, then for every $h\in\Hc_r$,
\beql{eq:center-Phi-general}
\Phi_r(h)=\sum_{\la\vdash r}\chi_q^\la(h)s_\la(u-v).
\eeq
\epr
\bpf
By Lemma~\ref{lem:center-rank-idempotent}, we may specialize the \emph{numerical idempotent} $E_\La$ at $q=1$ in the expression for $\Phi_r(e_\La)$ while keeping the alphabet variables independent.
The tensor Hecke action then becomes the signed superpermutation representation of $\Sym_r$.
For a permutation $\sigma$ with cycle lengths $\ell_1,\ldots,\ell_t$,
\beql{eq:center-cycle-supertrace}
\str_{V^{\ot r}}\big(Z^{\ot r}P_{\rm s}(\sigma)\big)
=\prod_{a=1}^t p_{\ell_a}(u-v),
\eeq
where $p_k(u-v)=\sum_i u_i^k-\sum_jv_j^k$.
Indeed a diagonal tensor contribution from a cycle of length $\ell$ requires all entries in that cycle to be equal.
For an odd entry the permutation contributes the sign $(-1)^{\ell-1}$ and the tensor supertrace contributes $(-1)^\ell$; their product is $-1$, independently of $\ell$.
For an even entry, the sign is $+1$.  Disjoint cycles multiply, proving \eqref{eq:center-cycle-supertrace}.

Let $f^\la$ be the number of standard $\la$-tableaux and $\chi^\la$ the irreducible character of $\Sym_r$.  
The central idempotent of its $\la$-block is
\[
 e^{\rm cen}_\la=\frac{f^\la}{r!}\sum_{\sigma\in\Sym_r}\chi^\la(\sigma^{-1})\sigma.
\]
All diagonal primitive idempotents in that block have equal weighted supertrace: matrix-unit cyclicity identifies the trace of $e^\la_{\La}$ with that of $e^\la_{\Ga}$.
Taking $1/f^\la$ times the central-idempotent trace, and using
\eqref{eq:center-cycle-supertrace}, gives
\[
\Phi_r(e_\La)=\frac1{r!} \sum_{\sigma\in\Sym_r}\chi^\la(\sigma^{-1})\prod_{a=1}^t p_{\ell_a(\sigma)}(u-v).
\]
Thus, we get \eqref{eq:center-Phi-Schur} by \eqref{eq:super-Frobenius}.

Finally, choose Hecke matrix units $e^\la_{\Ga\La}$ over $\CC(q)$.
Their weighted traces vanish when $\Ga\ne \La$, because $e^\la_{\Ga\Ga}e^\la_{\Ga\La}=e^\la_{\Ga\La}$ while $e^\la_{\Ga\La}e^\la_{\Ga\Ga}=0$, and $\Phi_r$ is cyclic.
Every diagonal unit has trace $s_\la(u-v)$.
If $h=\sum_{\la,\Ga,\La}h^\la_{\Ga\La}e^\la_{\Ga\La}$, then $\chi_q^\la(h)=\sum_\La h^\la_{\La\La}$, proving \eqref{eq:center-Phi-general}.
\epf

\bco\label{cor:center-homogeneous-HC}
If $\La$ has shape $\la\vdash r$, then
\beql{eq:center-homogeneous-HC}
\mathrm{pr}(C_\La)=s_\la(u-v),\qquad E_\La\ne0\ \Longleftrightarrow\ \la\in H_{M,N}.
\eeq
For arbitrary $h\in\Hc_r$,
\beql{eq:center-general-character}
\mathrm{pr}\,\str_{q,1,\ldots,r}\left(L_{\bar1}\cdots L_{\bar r}\varrho_r(h)\right)
=\sum_{\la\vdash r}\chi_q^\la(h)s_\la(u-v).
\eeq
\eco
\bpf
The first equality and the general character identity follow from Lemma~\ref{lem:center-Cartan-content} and Proposition~\ref{prop:center-weighted-Frob}. 
If $\la$ lies outside the hook, $s_\la(u-v)=0$.
Each content term in \eqref{eq:center-image-content} is a distinct monomial with
coefficient $(-1)^{\bar{\boldsymbol a}}$ times the nonnegative integer rank of the restricted idempotent.
Vanishing of the whole polynomial forces each such rank to vanish; hence $E_\La=0$.
The converse holds because a hook Schur polynomial is nonzero.
\epf

\subsection{The classical dimension bound and the whole center}\label{subsec:center-exhaustion}
We now combine the explicit central characters with the integral specialization. 
Write $F_s\Ucirc$ for the word filtration in the entries of $L$, equivalently in the entries of $X$; the two sets of generators are related by the affine change \eqref{eq:center-modified-X} and therefore define the same filtration over $\CC(q)$.

For each $s\ge0$ define a homomorphism of finite free $\Ac$-modules
\beql{eq:center-Ds}
\begin{aligned}
\Dc_s:F_s\Bc_{\Ac}&\longrightarrow\bigoplus_{i,j=1}^{M+N}F_{s+1}\Bc_{\Ac},\\
b&\longmapsto\big(bX_{ij}-X_{ij}b\big)_{i,j}.
\end{aligned}
\eeq
Since the $X_{ij}$ generate $\Bc_{\Ac}$, an element is in the kernel precisely when it commutes with every generator.
After extension of scalars, Theorem~\ref{thm:center-local-PBW} gives
\beql{eq:center-generic-kernel}
\ker(\Dc_s\ot_{\Ac}\CC(q))=Z(\Ucirc)\cap F_s\Ucirc
\eeq
and
\beql{eq:center-special-kernel}
\ker(\Dc_s\ot_{\Ac}\CC)=Z(U(\gl_{M|N}))\cap F_sU(\gl_{M|N}).
\eeq
The second identification uses isomorphism \eqref{eq:center-filter-specialization}: the images of the generators $X_{ij}$ are nonzero signed multiples of the standard $E_{ji}$ and thus generate the entire classical enveloping algebra.

\ble\label{lem:center-kernel-specialization}
The dimensions of the two kernels satisfy
\beql{eq:center-kernel-inequality}
\dim_{\CC(q)}\ker(\Dc_s\ot_{\Ac}\CC(q))\leq\dim_{\CC}\ker(\Dc_s\ot_{\Ac}\CC).
\eeq
\ele
\bpf
Choose finite ordered-PBW $\Ac$-bases of the source and target of $\Dc_s$.  The map is represented by a rectangular matrix $D_s(q)$ with entries in $\Ac$. 
Let $d_s$ be the rank of the source as a free $\Ac$-module; both of its fibers have dimension
$d_s$. 
Put
\[
r_0=\operatorname{rank}_{\CC}D_s(1).
\]
If $r_0>0$, the definition of matrix rank supplies an $r_0\times r_0$ minor whose determinant at $q=1$ is nonzero.
Taking the same rows and columns in $D_s(q)$ gives a determinant $\Delta(q)\in\Ac$ satisfying $\Delta(1)\ne0$.
The unique maximal ideal of $\Ac$ is $(q-1)\Ac$; hence $\Delta(q)$ is a unit of $\Ac$, in particular nonzero in $\CC(q)$.
Therefore
\[
\operatorname{rank}_{\CC(q)}D_s(q)\geq\operatorname{rank}_{\CC}D_s(1)=r_0.
\]
The assertion is automatic if $r_0=0$.
Rank--nullity now gives
\[
\begin{aligned}
\dim_{\CC(q)}\ker(D_s(q))&=d_s-\operatorname{rank}_{\CC(q)}D_s(q)\\
&\leq d_s-\operatorname{rank}_{\CC}D_s(1)=\dim_{\CC}\ker(D_s(1)),
\end{aligned}
\]
which proves the claim.
\epf

Combining the preceding lemma with the classical Harish--Chandra bound \eqref{eq:center-classical-bound} gives
\beq
\begin{aligned}
\dim_{\CC(q)}\big(Z(\Ucirc)\cap F_s\Ucirc\big)&\leq \dim_{\CC}\big(Z(U(\gl_{M|N}))\cap F_sU(\gl_{M|N})\big)\\
&\leq h_{M,N}^s.
\end{aligned}
\label{eq:center-upper-bound}
\eeq
We now show that the inequality is an equality, using the explicit homogeneous central elements constructed above.
\bth\label{thm:center-main}
For generic $q$, the ordinary center of the reflection subalgebra $\Ucirc$ is purely even, and
\beql{eq:center-intersection}
 Z(\Ucirc)=\Ucirc\cap Z(\U_q(\gl_{M|N})).
\eeq
For each $\la\in H_{M,N}$ choose any standard tableau $\La$ of shape $\la$, and write $C_\la=C_\La$.  Then $C_\la$ is independent of this choice and
\beql{eq:center-basis}
Z(\Ucirc)\cap F_s\Ucirc=\bigoplus_{\substack{\la\in H_{M,N}\\|\la|\le s}}\CC(q) C_\la,\qquad
Z(\Ucirc)=\bigoplus_{\la\in H_{M,N}}\CC(q) C_\la.
\eeq
The PBW Cartan projection restricts to an algebra isomorphism
\beql{eq:center-HC-isomorphism}
\HC=\mathrm{pr}|_{Z(\Ucirc)}:Z(\Ucirc)\overset{\sim}{\longrightarrow}\Lambda_{M|N},
\qquad \HC(C_\la)=s_\la(u-v),
\eeq
with $u_i,v_j$ specialized as in \eqref{eq:center-uv}.
\eth
\bpf
For each hook partition $\la$ the element $C_\La$ is even and central in the entire quantum superalgebra by Lemma~\ref{lem:center-monodromy}. 
Since $C_\La$ is the quantum supertrace of an ordered product of $|\la|$ matrices $L$, it
belongs to $F_{|\la|}\Ucirc$. 
By Corollary~\ref{cor:center-homogeneous-HC} its PBW projection is the nonzero polynomial $s_\la(u-v)$. 
By Proposition~\ref{prop:center-classical-input}, the images $\{s_\la(u-v):\la\in H_{M,N},\ |\la|\le s\}$ are linearly independent.
Therefore at least $h_{M,N}^s$ independent elements lie in $Z(\Ucirc)\cap F_s\Ucirc$.
The upper bound \eqref{eq:center-upper-bound} shows that there can be no further elements in that filtered center, and proves the first equality of \eqref{eq:center-basis}.

The word filtration exhausts $\Ucirc$, so taking the union over $s$ yields the second equality.
In particular every element of $Z(\Ucirc)$ is even and is an ambient central element; this proves the inclusion from left to right in \eqref{eq:center-intersection}. 
The reverse inclusion is immediate from the definition of the center of a subalgebra.
As recalled in Section~\ref{subsec:hc-pbw-prelim}, the PBW Cartan projection is multiplicative on $Z(\U_q(\gl_{M|N}))$.
Since we have now proved that all central elements of $\Ucirc$ are central in $\U_q(\gl_{M|N})$, $\HC$ is an algebra homomorphism.
On the basis \eqref{eq:center-basis} it takes the values $s_\la(u-v)$, which form a basis of $\Lambda_{M|N}$; therefore it is bijective.

Finally, if $\La,\La'$ have the same hook shape, their central elements $C_\La,C_{\La'}$ have identical Cartan projections by \eqref{eq:center-homogeneous-HC}. 
Injectivity of the isomorphism $\HC$ implies $C_\La=C_{\La'}$.
If the shape lies outside the hook, $E_\La=0$ by Corollary~\ref{cor:center-homogeneous-HC}.
Thus the notation $C_\la$ is unambiguous and completes the proof.
\epf

\bre\label{rem:center-ambient-comparison}
The description \eqref{eq:center-HC-isomorphism} is a theorem about the reflection subalgebra $\Ucirc$,
not a description of the full center of $\U_q(\gl_{M|N})$. The latter permits negative Cartan powers.  In particular, $z_0=\prod_{i=1}^M l_{ii}^+\prod_{j=1}^N (l_{M+j,M+j}^+)^{-1}$ is a central group-like element of $\U_q(\gl_{M|N})$. 
One can calculate
\[
z_0^2=q^{(M-N)(M-N-1)}\frac{u_1\cdots u_M}{v_1\cdots v_N}.
\]
When $N>0$, this is not a polynomial in $u,v$. 
Consequently $z_0\notin\Ucirc$: otherwise its square would belong to $Z(\Ucirc)$ and its $\HC$ image would violate \eqref{eq:center-HC-isomorphism}.
\ere

\section{Super $q$-immanants}\label{sec:q-immanants}
For the Hecke algebra $\Hc_m$ we use the representation $\varrho_m$ of \eqref{Hecke-action}, with $m$ the number of boxes of the partition under consideration. Under this representation,
\beql{eq:Hecke-baxter-action}
T_k(x,y)\longmapsto \Rcheck_k(x/y),\qquad
\Rcheck_k(u)=\Rcheck_k+\frac{q-q^{-1}}{u-1}.
\eeq

Let $\La$ be a standard tableau of shape $\la\vdash m$.
In what follows we assume that $E_\La\neq0$; if $E_\La=0$, all the formulas involving $E_\La$ below are identically zero.

We use the braided copies $L_{\bar a}$ defined in Section~\ref{subsec:mod-ref-gen}.

Define the polynomial $S_\La(z)$ by
\beql{eq:SLa-definition}
\begin{aligned}
S_\La(z)=\str_{1,\ldots,m}&(L^+_1+zq^{2c_1(\La)}L^-_1)\cdots(L^+_m+zq^{2c_m(\La)}L^-_m)\\
&\times (L^-_m)^{-1}\cdots(L^-_1)^{-1}D_1\cdots D_mE_\La.
\end{aligned}
\eeq

\bpr\label{prop:equivalent-form}
We have
\beql{eq:SLa-equivalent}
S_\La(z)=\str_{q;1,\ldots,m}(L_{\bar1}+zq^{2c_1(\La)})\cdots(L_{\bar m}+zq^{2c_m(\La)})E_\La.
\eeq
\epr
\bpf
Substitute $t_a=zq^{2c_a(\La)}$ into the product identity \eqref{eq:center-shifted-product} and take the quantum supertrace.
\epf

Let $L_q(\eta|\theta)$ be a highest-weight $\U_q(\gl_{M|N})$-module with highest weight
$(\eta_1,\ldots,\eta_M\mid\theta_1,\ldots,\theta_N)$.
A highest-weight vector $\xi$ satisfies
\[
l^-_{ij}\xi=0\quad(i>j),\qquad l^+_{ii}\xi=q^{\eta_i}\xi\quad(1\leq i\leq M),\qquad
l^+_{M+j,M+j}\xi=q^{-\theta_j}\xi\quad(1\leq j\leq N).
\]
Every element of $Z(\U_q(\gl_{M|N}))$ acts on this module by a scalar obtained by evaluating its Harish--Chandra image at the specified highest weight.

A tableau $T$ of shape $\mu$ with entries in $\{1,\ldots,M+N\}$ is called $(M,N)$-semistandard if its entries weakly increase along each row and down each column, while the entries $1,\ldots,M$ increase strictly down each column and the entries $M+1,\ldots,M+N$ increase strictly along each row.

Let
\[
x=(x_1,\ldots,x_M),\qquad
y=(y_1,\ldots,y_N),
\]
and let $a=(a_k)_{k\in\BZ}$ be a doubly infinite sequence of parameters.

For a partition $\mu$ contained in the $(M,N)$-hook, define the factorial supersymmetric Schur polynomial by
\beq
s_\mu^{M|N}(x/y\mid a)=\sum_{\mathcal T}\prod_{\alpha\in\mu}w_{\mathcal T}(\alpha),
\eeq
where the sum runs over all $(M,N)$-semistandard tableaux $\mathcal T$ of shape $\mu$, and
\beq
w_{\mathcal T}(\alpha)=
\begin{cases}
x_i+a_{i+c(\alpha)},&\mathcal T(\alpha)=i,\quad 1\leq i\leq M,\\
y_j-a_{M-j+1+c(\alpha)},&\mathcal T(\alpha)=M+j,\quad 1\leq j\leq N.
\end{cases}
\eeq

\bth
\label{thm:main}
Let $\La$ be a standard tableau of shape $\la\vdash m$.
\begin{enumerate}[label=\textup{\arabic*.}]
\item All coefficients of the polynomial $S_\La(z)$ belong to the center
of $\U_q(\gl_{M|N})$.

\item The polynomial $S_\La(z)$ does not depend on the choice of the
standard $\la$-tableau $\La$.  Hence it will be denoted by $S_\la(z)$.

\item
For a partition $\mu$ in the $(M,N)$-hook and a highest-weight module $L_q(\eta\mid\theta)$, the eigenvalue of $S_\mu(z)$ is $s_\mu^{M|N}(x/y\mid a)$, where
\[
x_i=q^{2(\eta_i-i+1)},\quad y_j=-q^{-2(\theta_j+M-j)},\quad a_r=zq^{2(r-M+N)}.
\]
for $1\leq i\leq M,1\leq j\leq N$ and $r\in\BZ$.
\end{enumerate}
\eth

\bpf
Part~1 follows at once from Lemma~\ref{lem:center-monodromy} by taking $h=e_\La$ and $t_a=zq^{2c_a(\La)}$ in \eqref{eq:SLa-equivalent}.

We now prove Part~2.  It is enough to consider two standard tableaux related by an adjacent transposition. Let
\[
\La'=\si_k\La,\qquad d_k=c_{k+1}(\La)-c_k(\La),
\]
and assume that $\La'$ is standard.
The RTT relation \eqref{ev RTT} gives
\beql{eq:tableau-intertwiner}
L_{\La'}(z)\Rcheck_k(q^{-2d_k})
=\Rcheck_k(q^{-2d_k})L_\La(z).
\eeq
On the other hand, Lemma~\ref{lem:et} implies
\beql{eq:idempotent-intertwiner}
\Rcheck_k(q^{-2d_k})E_\La
=E_{\La'}\Rcheck_k(q^{2d_k}).
\eeq
Both $\Rcheck_k(q^{-2d_k})$ and $\Rcheck_k(q^{2d_k})$ are invertible.
Since $D_1\cdots D_m$ commutes with $\Rcheck_k(u)$, equations \eqref{eq:tableau-intertwiner} and \eqref{eq:idempotent-intertwiner} reduce the two expressions for $S_\La(z)$ and $S_{\La'}(z)$ to conjugate operators under the supertrace. 
The cyclicity of the supertrace for even operators therefore yields
\[
S_{\La'}(z)=S_\La(z).
\]
Any two standard tableaux of the same shape can be connected by a sequence of admissible adjacent transpositions. 
Hence $S_\La(z)$ depends only on $\la$.

Finally, we prove Part~(3). By Part~(2), the polynomial $S_\La(z)$ depends only on the shape $\mu$, so we write it as $S_\mu(z)$.

By Theorem 5.2 in \cite{JLZ 2026}, we have
\beql{eq:HC-reverse-tableau}
\chi(S_\mu(z))=\sum_{\mathcal T}(-1)^{\bar{\mathcal T}}\prod_{\alpha\in\mu}X_{\mathcal T(\alpha)}
\bigl(zq^{2c(\alpha)}\bigr),
\eeq
where the sum is over all reverse $(M,N)$-semistandard tableaux of shape $\mu$, and
\[
\bar{\mathcal T}=\sum_{\alpha\in\mu}\overline{\mathcal T(\alpha)}\pmod2.
\]
Here the entries weakly decrease along rows and columns, with strict decrease down columns for even entries and along rows for odd entries.

It remains to calculate the functions $X_i(u)$. Recall that
\[
D=\diag(d_1,\ldots,d_{M+N}),\qquad d_i=
\begin{cases}
q^{-2i+2},&1\leq i\leq M,\\
q^{-2M+2j},&i=M+j,\quad 1\leq j\leq N.
\end{cases}
\]
By the highest-weight conditions,
\[
l^-_{ij}\xi=0\quad(i>j),
\]
and
\[
l^+_{ii}\xi=q^{\eta_i}\xi,\qquad l^+_{M+j,M+j}\xi=q^{-\theta_j}\xi.
\]
The diagonal contribution of
\[
(L^++uL^-)(L^-)^{-1}D
\]
on the highest-weight vector is therefore
\[
\begin{aligned}
X_i(u)&=q^{-2i+2}(q^{2\eta_i}+u),&&1\leq i\leq M,\\
X_{M+j}(u)&=q^{-2M+2j}(q^{-2\theta_j}+u),&&1\leq j\leq N.
\end{aligned}
\]

For $1\leq i\leq M$ and a box $\alpha$, we have
\[
\begin{aligned}
X_i\bigl(zq^{2c(\alpha)}\bigr)
&=q^{-2i+2}
\left(q^{2\eta_i}+zq^{2c(\alpha)}\right)\\
&=x_i+zq^{2(c(\alpha)-i+1)}\\
&=x_i+a_{M-N+1-i+c(\alpha)}.
\end{aligned}
\]
Similarly, for $1\leq j\leq N$, the sign
from the supertrace gives
\[
\begin{aligned}
-X_{M+j}\bigl(zq^{2c(\alpha)}\bigr)
&=-q^{-2M+2j}
\left(q^{-2\theta_j}+zq^{2c(\alpha)}\right)\\
&=y_j-zq^{2(c(\alpha)-M+j)}\\
&=y_j-a_{j-N+c(\alpha)}.
\end{aligned}
\]

Consequently, the sign $(-1)^{\bar{\mathcal T}}$
can be absorbed into the factors corresponding
to entries $M+j$, $1\leq j\leq N$, and we obtain
\[
\chi(S_\mu(z))=\sum_{\mathcal T}\prod_{\substack{\alpha\in\mu\\\mathcal T(\alpha)=i}}
\left(x_i+a_{M-N+1-i+c(\alpha)}\right)\prod_{\substack{\alpha\in\mu\\
\mathcal T(\alpha)=M+j}}\left(y_j-a_{j-N+c(\alpha)}\right).
\]
By the equivalent reverse-tableau formula
for factorial supersymmetric Schur polynomials,
which follows from the skew-factorial
decomposition in \cite{Molev 1998},
the right-hand side equals
\[
s_\mu^{M|N}(x/y\mid a).
\]
This completes the proof of Part~(3).
\epf

\bco\label{cor:super-JLM-homogeneous}
For a hook partition $\la$, let $C_\la=S_\la(0)$.
The elements $\{C_\la:\la\in H_{M,N}\}$ form a $\CC(q)$-basis of $Z(\Ucirc)$, and
\[
\mathrm{pr}(C_\la)=s_\la(u-v),\qquad u_i=d_i(l_{ii}^+)^2,\quad v_j=d_{M+j}(l_{M+j,M+j}^+)^2.
\]
For $\la\notin H_{M,N}$ we have $C_\la=0$.
\eco
\bpf
By Proposition~\ref{prop:equivalent-form}, $S_\La(0)=C_\La$. The assertions follow from Corollary~\ref{cor:center-homogeneous-HC} and Theorem~\ref{thm:center-main}.
\epf

\section{Higher super Capelli identities}\label{sec:super-capelli}
We extend the braided Weyl algebra construction of \cite[Section~4]{JLM 2025} to the standard super $R$-matrix.
The matrix Capelli identity is covered by the general Hecke-symmetry framework \cite{GPS 2022,GSZ 2025,Zaitsev 2024}; the purpose here is to identify its source in the reflection algebra $\Ucirc$ and keep the normalization of the modified generators in Section~\ref{subsec:mod-ref-gen}.

\subsection{The braided Weyl superalgebra}
Let $\mathcal W_{M|N}$ be the superalgebra generated by the entries
of two even total matrices $\mathsf M=(m_{ij})$ and
$\partial=(\partial_{ij})$, with
$|m_{ij}|=|\partial_{ij}|=\bar i+\bar j$, subject to
\beq
\begin{aligned}
\Rcheck\mathsf M_1\Rcheck\mathsf M_1
 &=\mathsf M_1\Rcheck\mathsf M_1\Rcheck,\\
\Rcheck^{-1}\partial_1\Rcheck^{-1}\partial_1
 &=\partial_1\Rcheck^{-1}\partial_1\Rcheck^{-1},\\
\partial_1\Rcheck\mathsf M_1\Rcheck
 &=\Rcheck\mathsf M_1\Rcheck^{-1}\partial_1+\Rcheck.
\end{aligned}
\label{eq:super-Weyl-relations}
\eeq
Here the products are taken in the graded tensor algebra. We
use the braided copies
\[
A_{\bar1}=A_1,\qquad
A_{\overline{a+1}}=\Rcheck_aA_{\bar a}\Rcheck_a^{-1},
\qquad A=\mathsf M,\partial,K.
\]

\bpr\label{prop:super-Weyl-realization}
Put $K=\mathsf M\partial$. Then $K$ satisfies
\beq
\Rcheck K_1\Rcheck K_1-K_1\Rcheck K_1\Rcheck
 =\Rcheck K_1-K_1\Rcheck.
\label{eq:super-Weyl-modified-RE}
\eeq
Consequently, there is a superalgebra homomorphism
\beq
\rho:\Ucirc\longrightarrow\mathcal W_{M|N},\qquad
\rho(X)=K,\qquad
\rho(L)=I-(q-q^{-1})\mathsf M\partial.
\label{eq:super-Weyl-map}
\eeq
\epr
\bpf
The mixed relation in \eqref{eq:super-Weyl-relations} gives
\[
\partial_1\Rcheck\mathsf M_1
=\Rcheck\mathsf M_1\Rcheck^{-1}
\partial_1\Rcheck^{-1}+I.
\]
Using the other two relations yields
\[
\begin{aligned}
\Rcheck K_1\Rcheck K_1
 &=\mathsf M_1\Rcheck\mathsf M_1
\partial_1\Rcheck^{-1}\partial_1+\Rcheck K_1,\\
K_1\Rcheck K_1\Rcheck
 &=\mathsf M_1\Rcheck\mathsf M_1
\partial_1\Rcheck^{-1}\partial_1+K_1\Rcheck.
\end{aligned}
\]
Subtracting proves \eqref{eq:super-Weyl-modified-RE}.
By Theorem~\ref{thm:center-local-PBW}, the presentation
\eqref{eq:center-modified-RE} identifies the abstract modified
reflection algebra over $\CC(q)$ with $\Ucirc$. Thus it defines
\eqref{eq:super-Weyl-map}.
\epf

\subsection{The Capelli normalization and matrix identity}
For $c\in\BZ$ set
\beq
\kappa(c)=\frac{1-q^{-2c}}{q-q^{-1}},
\qquad
\kappa(c)\big|_{q=1}=c.
\label{eq:super-capelli-kappa}
\eeq
For a standard tableau $\La$ of shape $\la\vdash r$, put
\beq
\mathcal C_\La
=\str_{q;1,\ldots,r}
\left(\prod_{a=1}^r
\bigl(X_{\bar a}-\kappa(c_a(\La))\bigr)E_\La\right),
\qquad
X_{\bar a}=\frac{I-L_{\bar a}}{q-q^{-1}}.
\label{eq:super-capelli-elements}
\eeq
The general centrality assertion of Lemma~\ref{lem:center-monodromy}
shows that $\mathcal C_\La\in Z(\Ucirc)$.
For the content convention used in the present paper, the
immanants $S_\La(z)$ of \eqref{eq:SLa-equivalent} have shifts
$q^{2c_a(\La)}$. The Capelli normalization instead involves
$q^{-2c_a(\La)}$; explicitly,
\beq
\mathcal C_\La
=(-1)^r(q-q^{-1})^{-r}
\str_{q;1,\ldots,r}
\left(\prod_{a=1}^r(L_{\bar a}-q^{-2c_a(\La)})E_\La\right).
\label{eq:super-capelli-parameter-comparison}
\eeq
Thus the negative-content Capelli family in this section is kept
separate from the positive-content family $S_\La(z)$ of Section~6.
Its matrix identity below does not require an identification
with $S_\la(-1)$ or a new tableau-independence assertion.

Let $J_1=1$, $J_{a+1}=T_aJ_aT_a$ be the Hecke Jucys--Murphy
elements introduced in Section~\ref{subsec:center-idempotents-prelim},
and set $\mathsf J_a=\varrho_r(J_a)$. Then
\[
\mathsf J_r^{-1}E_\La=q^{-2c_r(\La)}E_\La.
\]
Write
\[
\mathsf N_r=\mathsf M_{\bar1}\cdots\mathsf M_{\bar r}
\partial_{\bar r}\cdots\partial_{\bar1},\qquad
\mathsf N_0=I.
\]
The following contraction identity is the standard Hecke-type
normal-ordering relation; see \cite[Section~4]{JLM 2025}
and \cite{GSZ 2025}.
\ble\label{lem:super-capelli-contraction}
For $r\ge2$ one has
\beq
\partial_{\overline{r-1}}\cdots\partial_{\bar1}K_{\bar r}
=K_{\bar r}\partial_{\overline{r-1}}\cdots\partial_{\bar1}
\mathsf J_r^{-1}
+\partial_{\overline{r-1}}\cdots\partial_{\bar1}
\frac{I-\mathsf J_r^{-1}}{q-q^{-1}}.
\label{eq:super-capelli-contraction}
\eeq
Consequently,
\beq
\mathsf N_{r-1}K_{\bar r}
=\mathsf N_r\mathsf J_r^{-1}
+\mathsf N_{r-1}\frac{I-\mathsf J_r^{-1}}{q-q^{-1}}.
\label{eq:super-capelli-recursion}
\eeq
\ele
\bpf
The mixed relation in \eqref{eq:super-Weyl-relations}
gives, for two adjacent auxiliary positions,
\[
\partial_{\bar1}K_{\bar2}
=K_{\bar2}\partial_{\bar1}\Rcheck_1^{-2}
+\partial_{\bar1}\Rcheck_1^{-1}.
\]
Transporting this identity along the braid words and moving
$K_{\bar r}$ through the derivative factors gives a leading
coefficient $\mathsf J_r^{-1}$.
The remaining coefficients are a telescoping sum:
if $A_a$ are the braid factors of
\cite[proof of Theorem~4.1]{JLM 2025}, their contraction partners
are $(I-A_a)/(q-q^{-1})$, and hence
\[
\sum_{a=1}^{r-1}(I-A_a)A_{a-1}\cdots A_1
=I-A_{r-1}\cdots A_1=I-\mathsf J_r^{-1}.
\]
This proves \eqref{eq:super-capelli-contraction}; multiplying
by $\mathsf M_{\bar1}\cdots\mathsf M_{\overline{r-1}}$
gives \eqref{eq:super-capelli-recursion}.
\epf

\bth\label{thm:super-capelli}
For every standard tableau $\La$ of shape $\la\vdash r$,
\beq\label{eq:super-capelli-matrix}
\left(\prod_{a=1}^r
(K_{\bar a}-\kappa(c_a(\La)))\right)E_\La
=\left(\prod_{\alpha\in\la}q^{-2c(\alpha)}\right)
\mathsf M_{\bar1}\cdots\mathsf M_{\bar r}
\partial_{\bar r}\cdots\partial_{\bar1}E_\La.
\eeq
Consequently, the image of the central element
\eqref{eq:super-capelli-elements} under $\rho$ is
\beq
\rho(\mathcal C_\La)=
\left(\prod_{\alpha\in\la}q^{-2c(\alpha)}\right)
\str_{q;1,\ldots,r}\left(
\mathsf M_{\bar1}\cdots\mathsf M_{\bar r}
\partial_{\bar r}\cdots\partial_{\bar1}E_\La\right).
\label{eq:super-capelli-scalar}
\eeq
\eth
\bpf
For $r=1$ the assertion is $K=\mathsf M\partial$.
Let $\La^-$ be obtained by removing the box containing $r$.
The primitive idempotents satisfy $E_{\La^-}E_\La=E_\La$,
and $K_{\bar r}$ commutes with the Hecke action in the first
$r-1$ positions. By induction, the left side of
\eqref{eq:super-capelli-matrix} is
\[
\left(\prod_{\alpha\in\mathrm{sh}(\La^-)}
q^{-2c(\alpha)}\right)
\mathsf N_{r-1}\bigl(K_{\bar r}-\kappa(c_r(\La))\bigr)E_\La.
\]
Using \eqref{eq:super-capelli-recursion} and
$\mathsf J_r^{-1}E_\La=q^{-2c_r(\La)}E_\La$,
we obtain
\[
\begin{aligned}
&\mathsf N_{r-1}\bigl(K_{\bar r}-\kappa(c_r(\La))\bigr)E_\La\\
&\quad=q^{-2c_r(\La)}\mathsf N_rE_\La
+\mathsf N_{r-1}
\left(\frac{1-q^{-2c_r(\La)}}{q-q^{-1}}
-\kappa(c_r(\La))\right)E_\La\\
&\quad=q^{-2c_r(\La)}\mathsf N_rE_\La.
\end{aligned}
\]
The scalar factors multiply to the one in
\eqref{eq:super-capelli-matrix}. Taking the quantum supertrace
proves \eqref{eq:super-capelli-scalar}.
\epf

\subsection{Classical specialization}
The relations \eqref{eq:super-Weyl-relations} have a regular
specialization at $q=1$. There $\Rcheck$ becomes the graded flip,
so the coordinate and derivative subalgebras become
supercommutative, while the mixed relation becomes
\[
\partial_1\mathsf M_2=\mathsf M_2\partial_1+P.
\]
Thus $\mathcal W_{M|N}$ specializes to the Weyl superalgebra
of polynomial differential operators. Moreover,
$\kappa(c)\big|_{q=1}=c$ and
$E_\La$ specializes to the usual Young idempotent. Consequently,
\eqref{eq:super-capelli-matrix} is a $q$-deformation of the
higher super Capelli identity associated with $\la$.
The source algebra specializes through
$X_{ij}\mapsto-(-1)^{\bar j}E_{ji}$, as in
Lemma~\ref{lem:center-special-fiber}. Thus the classical
identity is expressed in the same matrix-unit convention as the
classical super Capelli formulas of \cite{Molev 1998}.
The operators $\rho(\mathcal C_\La)$ commute with $\rho(\Ucirc)$;
centrality in the whole braided Weyl superalgebra is not asserted.

\end{document}